\documentclass[12pt]{amsart}

\usepackage{amsmath}
\usepackage{amscd}
\usepackage{amssymb}
\usepackage{latexsym}
\usepackage{stmaryrd} 
\usepackage{mathbbol} 

\input xy
\xyoption{all} \CompileMatrices

\newtheorem{theorem}{Theorem}[section]
\newtheorem*{theorem*}{Theorem}
\newtheorem{lemma}[theorem]{Lemma}
\newtheorem*{lemma*}{Lemma}
\newtheorem{corollary}[theorem]{Corollary}
\newtheorem{proposition}[theorem]{Proposition}

\newtheorem{remark}[theorem]{Remark}
\newtheorem{definition}[theorem]{Definition}

\newcommand{\bgl}{\begin{equation}} 
\newcommand{\egl}{\end{equation}}
\newcommand{\bgloz}{\begin{equation*}} 
\newcommand{\egloz}{\end{equation*}}
\newcommand{\bgln}{\begin{eqnarray}} 
\newcommand{\egln}{\end{eqnarray}}
\newcommand{\bglnoz}{\begin{eqnarray*}} 
\newcommand{\eglnoz}{\end{eqnarray*}}
\newcommand{\btheo}{\begin{theorem}}
\newcommand{\etheo}{\end{theorem}}
\newcommand{\btheooz}{\begin{theorem*}}
\newcommand{\etheooz}{\end{theorem*}}
\newcommand{\blemma}{\begin{lemma}}
\newcommand{\elemma}{\end{lemma}}
\newcommand{\blemmaoz}{\begin{lemma*}}
\newcommand{\elemmaoz}{\end{lemma*}}
\newcommand{\bproof}{\begin{proof}}
\newcommand{\eproof}{\end{proof}}
\newcommand{\bbew}{\begin{beweis}}
\newcommand{\ebew}{\end{beweis}}
\newcommand{\bremark}{\begin{remark}\em}
\newcommand{\eremark}{\end{remark}}
\newcommand{\bdefin}{\begin{definition}}
\newcommand{\edefin}{\end{definition}}
\newcommand{\bprop}{\begin{proposition}}
\newcommand{\eprop}{\end{proposition}}
\newcommand{\bcor}{\begin{corollary}}
\newcommand{\ecor}{\end{corollary}}
\newcommand{\bfa}{\begin{cases}} 
\newcommand{\efa}{\end{cases}}
\newcommand{\cA}{\mathcal A}

\newcommand{\cC}{\mathcal C}

\newcommand{\cF}{\mathcal F}

\newcommand{\cH}{\mathcal H}
\newcommand{\cI}{\mathcal I}

\newcommand{\cL}{\mathcal L}
\newcommand{\cM}{\mathcal M}

\newcommand{\cO}{\mathcal O}

\newcommand{\cQ}{\mathcal Q}

\newcommand{\cU}{\mathcal U}

\newcommand{\so}{\text{\tiny{$\mathcal O$}}}

\def\Az{\mathbb{A}}
\def\Cz{\mathbb{C}}

\def\Nz{\mathbb{N}}

\def\Qz{\mathbb{Q}}

\def\Zz{\mathbb{Z}}

\def\1z{\mathbb{1}}
\newcommand{\fA}{\mathfrak A}

\newcommand{\fD}{\mathfrak D}

\newcommand{\fI}{\mathfrak I}

\newcommand{\an}[1]{``#1''} 
\newcommand{\ti}{\tilde}

\newcommand{\ri}{\rightarrow}
\newcommand{\lori}{\longrightarrow}

\newcommand{\rarr}{\rightarrow}

\newcommand{\ma}{\mapsto} 
\newcommand{\loma}{\longmapsto} 
\newcommand{\lomafr}{\longmapsfrom} 
\newcommand\onto{\twoheadrightarrow} 
\newcommand\into{\hookrightarrow} 
\newcommand{\Rarr}{\Rightarrow} 
\newcommand{\LRarr}{\Leftrightarrow} 

\def\SEMI{\mbox{$\times\kern-2pt\vrule height5pt width.6pt \kern3pt $}}

\newcommand{\halb}{\tfrac{1}{2}}

\newcommand{\End}{{\rm End}\,}
\newcommand{\Aut}{{\rm Aut}\,}

\newcommand{\img}{{\rm Im\,}}

\newcommand{\Spec}{{\rm Spec\,}} 

\newcommand{\alg}{{\rm alg}}
\newcommand{\de}{{\rm deg\,}} 
\newcommand{\ntei}{\nmid} 

\newcommand{\Ad}{{\rm Ad\,}}

\newcommand{\reg}{^\times} 
\newcommand{\pos}{_{>0}} 
\newcommand{\co}{^{\complement}} 
\newcommand{\ev}{\operatorname{ev}} 
\newcommand{\lspan}{{\rm span}} 
\newcommand{\clspan}{\overline{\lspan}} 
\newcommand{\Iso}{{\rm Isom}} 
\newcommand{\norm}[1]{\left\|#1\right\|} 
\newcommand{\defeq}{\mathrel{:=}} 

\newcommand{\dop}{\text{: }} 
\newcommand{\falls}{\text{ if }} 
\newcommand{\sonst}{\text{ else}} 
\newcommand{\fa}{\text{ for all }} 
\newcommand{\ilim}{\varinjlim} 
\newcommand{\plim}{\varprojlim} 
\newcommand{\e}[1]{e_{\left[#1\right]}} 
\newcommand{\E}[1]{E_{\left[#1\right]}} 
\newcommand{\chf}[1]{\1z_{\left[#1\right]}} 
\newcommand{\card}[1]{\# \left[#1\right]} 
\newcommand{\dotcup}{\ensuremath{\mathaccent\cdot\cup}} 
\newcommand{\rte}{\overset{e}{\rtimes}} 
\newcommand{\rtm}{\ti{\rtimes}} 
\newcommand{\lge}{\left\{} 
\newcommand{\rge}{\right\}} 
\newcommand{\lru}{\left(} 
\newcommand{\rru}{\right)} 
\newcommand{\leck}{\left[} 
\newcommand{\reck}{\right]} 
\newcommand{\rukl}[1]{\lru #1 \rru} 
\newcommand{\eckl}[1]{\leck #1 \reck} 
\newcommand{\gekl}[1]{\lge #1 \rge} 
\newcommand{\menge}[2]{\gekl{ #1 \dop #2 }} 
\newcommand{\ql}[2]{\rukl{#1} / #2} 
\newcommand{\qr}[2]{#1 / \rukl{#2}} 
\begin{document}

\title{Ring C*-algebras}

\author{Xin Li}


\thanks{\scriptsize{Research supported by the Deutsche Forschungsgemeinschaft and the Deutsche Telekom Stiftung.}}

\thanks{\scriptsize{This work has been done in the context of the author's PhD project at the University of Muenster.}}

\begin{abstract}
We associate reduced and full C*-algebras to arbitrary rings and study the inner structure of these ring C*-algebras. As a result, we obtain conditions for them to be purely infinite and simple. We also discuss several examples. Originially, our motivation comes from algebraic number theory.
\end{abstract}

\maketitle

\section{Introduction}

This paper continues the work of \cite{CuLi} on C*-algebras associated to rings. The original motivation behind our investigations came from algebraic number theory. It was the work of Bost and Connes, \cite{BoCo}, which initiated investigations of links between operator algebras and number theory. The main result of \cite{BoCo} was to construct a C*-dynamical system whose thermodynamical behaviour, described in terms of KMS states, reveals a close relationship to classfield theory over the rational numbers. This discovery has led many authors to investigate dynamical systems with similar properties in more general situations (see \cite{CoMa}, Chapter 3 for a survey of the developments).

Most relevant for us is the construction of Cuntz in \cite{Cun1}. His approach differs from the other ones in two main points:

Cuntz's investigation is focused on the C*-algebra itself rather than on the dynamical system. But this time, the construction really uses the full ring structure, not only the multiplicative part (see the explanation at the end of Section \ref{Re}). The resulting C*-algebra is still closely related to number theory, but at the same time, it has an interesting structure on its own, which is investigated in \cite{Cun1}. 

This is exactly the point where our story of ring C*-algebras begins, because this construction of Cuntz is nothing else but what we call the ring C*-algebra of the integers.

The first step of generalization has then been done in \cite{CuLi}, where we considered integral domains with finite quotients. It is precisely this finiteness condition which allows a straightforward generalization of Cuntz's construction. We will recall the construction of \cite{CuLi} in Section \ref{Re}. Still, even though this first step of generalization was very natural, it only covers a rather small class of rings. Consequently, the main issue of the present paper is the following question:

Is it possible to extend the construction of \cite{CuLi} to arbitrary rings?

Actually, one encounters a similar situation as going from $\cO_n$ for finite $n$ to $\cO_{\infty}$. But in our situation, we have many more generators which have to be organized reasonably. The central idea is to add one additional piece of data, namely a certain set-theoretical algebra over our ring. We will explain our new construction in detail in Section \ref{R}. The resulting C*-algebras are called ring C*-algebras in analogy to the group case. Moreover, we show that this new construction generalizes the former one of \cite{CuLi} in a very satisfactory way.

As a second step, we investigate the inner structure of these ring C*-algebras, first in the general case (Section \ref{B}), then in the special situation of commutative rings (Section \ref{C}). As a main result, we obtain necessary and sufficient conditions for ring C*-algebras to be purely infinite and simple. This result yields a characterization in terms of generators and relations of the reduced ring C*-algebra, which is a priori given as concrete operators on a Hilbert space (see Corollary \ref{A-Ar}).

Finally, we discuss some typical examples which illustrate our theory and reveal certain connections to algebraic number theory (Section \ref{E}).

We also mention that the main idea which enters into our construction of ring C*-algebras can be used to extend and clarify the existing theory of crossed products by semigroups. This is explained in the appendix, where the basics of semigroup crossed products are recalled as well.

\section{The finite case}

Let us recall the construction of \cite{CuLi}. With this special case in mind, we will then motivate and develop the notion of ring C*-algebras in the general case.

\subsection{Review}

\label{Re}

As one can see in \cite{CuLi}, it turns out that the construction as well as the structural analysis of \cite{Cun1} only used two properties of the integers: Namely, that $\Zz$ is an integral domain and that $\Zz$ has finite quotients. So here is the first step of generalization:

Let $R$ be an integral domain, which means that $R$ is a commutative, unital ring without zero-divisors. Moreover, assume that $R$ has finite quotients, in the sense that for all nontrivial $ b \in R $, we have $\card{R / (b)} < \infty$. We will always view our rings as purely algebraic objects, that is to say that we consider the discrete case only. 

We associate two C*-algebras to $R$:

\bdefin
Consider unitaries $\menge{U^a}{a \in R}$ and isometries $\menge{S_b}{b \in R \reg}$ (here $R \reg$ is $R \setminus \gekl{0}$) on the Hilbert space $\ell^2(R)$ given by
\bglnoz
  && U^a \xi_r = \xi_{a+r} \\
  && S_b \xi_r = \xi_{b \cdot r} 
\eglnoz
where $\menge{\xi_r}{r \in R}$ is the canonical orthonormal basis of $\ell^2(R)$. 

Then, the reduced ring C*-algebra of $R$ is given by 
\bgloz
  \fA_r [R] \defeq C^* \rukl{ \menge{U^a, S_b}{a \in R, b \in R \reg} }.
\egloz
\edefin

Note that this is very much in the spirit of reduced group C*-algebras. Moreover, we have already used the assumption that $R$ does not have zero-divisors: Otherwise the formula defining $S_b$ would not give rise to an isometry (it would not even produce a bounded operator in general).

At this point, to motivate the definition of the full ring C*-algebra, we make the following observation: 

The range projection of $S_b$ is given by $S_b S_b^*(\xi_r) = \chf{(b)}(r) \xi_r$ where $\chf{(b)}$ denotes the characteristic function on $(b) \subseteq R$. This follows immediately from the definitions. Similarly, $U^a S_b S_b^* U^{-a}$ corresponds to $ \chf{a+(b)} $ in the sense that $U^a S_b S_b^* U^{-a}(\xi_r) = \chf{a+(b)}(r) \xi_r$. This leads to a very important relation reflecting the fact that 
\bgloz
  R = \dotcup_{a+(b) \in R/(b)} a+(b),
\egloz
namely
\bgloz
  \sum_{a+(b) \in R/(b)} U^a S_b S_b^* U^{-a} = 1.
\egloz
This relation incorporates the ideal structure of $R$ in a natural way, and at the same time, it contains valuable information on the range projections of the $S_b$.

Thus, we proceed as follows:

\bdefin

\label{fin}

The full ring C*-algebra of $R$, denoted by $\fA [R]$, is the universal C*-algebra generated by
\bgloz
  \text{unitaries } \menge{u^a}{a \in R}
\egloz 
\bgloz
  \text{and isometries } \menge{s_b}{b \in R \reg}
\egloz 
satisfying the relations
\bglnoz
  I. && u^a s_b u^c s_d = u^{a+bc} s_{bd} \fa a,c \in R; b,d \in R \reg \\
  II. && \sum_{a+(b) \in R/(b)} u^a s_b s_b^* u^{-a} = 1 \fa b \in R \reg.
\eglnoz
\edefin

Note that I. implies that $u^a s_b s_b^* u^{-a}$ only depends on the coset $a+(b)$. Moreover, we have heavily used the finiteness condition in Relation II. We will come back to this.

The notion of universal C*-algebras is explained in \cite{Bla}, II.8.3. We just mention that in our case, existence of $\fA [R]$ is guaranteed as all the generators have norm less or equal to $1$. But since we already have a nontrivial realization of these generators and relations on $\ell^2(R)$, namely in $\fA_r [R]$, it follows that the full ring C*-algebra cannot be trivial. Actually, the universal property of $\fA [R]$ yields an epimorphism $\fA [R] \lori \fA_r [R]$ sending generators to generators. Once again, this is in the spirit of full group C*-algebras. 

Finally, we mention some properties of these constructions. As we will prove more general versions later on, proofs will be omitted at the moment.

First of all, $ \fA [R] $ admits a crossed product description $\fA [R] \cong \fD [R] \rte P_R$ where $ \fD [R] $ is the C*-subalgebra of $ \fA [R] $ generated by all the projections $ u^a s_b s_b^* u^{-a} $ (for $a \in R, b \in R \reg$) and $P_R = R \rtimes R \reg$ is the $ax+b$-semigroup over $R$. $ \fD [R] $ is commutative and the action of $P_R$ is given by  
\bgloz
  P_R \ni (a,b) \loma \Ad(u^a s_b) \in \End(\fD [R]).
\egloz
The crossed product is taken in the sense of Section \ref{semicropro}. This description already allows some consequences: It implies that $\fA [R]$ lies in the nuclear UCT class (these notions are explained in \cite{Ror}, 2.4).

With considerably more work, it can be shown that $ \fA [R] $ is purely infinite and simple if and only if $R$ is not a field. This is a very strong property, for instance, it immediately implies that the canonical surjection from $ \fA [R] $ onto $ \fA_r [R] $ must be faithful. 

For more details, we refer to \cite{CuLi}, where the relationship to generalized Bost-Connes algebras is explained as well. Moreover, just note that in \cite{Cun1}, $\fA [\Zz]$ is denoted by $\cQ_{\Zz}$. Cuntz also studies $ \cQ_{\Nz} $ which coincides with the C*-subalgebra 
\bgloz
  C^* \rukl{ \menge{u^a, s_b}{a \in \Zz,b \in \Zz \pos} } \subseteq \fA [\Zz].
\egloz

Moreover, Bost and Connes studied $C^* \rukl{\fD [\Zz], \menge{s_b}{b \in \Zz \pos}}$. That is why we said that their C*-algebra does not use the whole ring structure.

The close relationship to algebraic number theory is due to the following facts:

1. The rings of integers of algebraic number fields (or, more generally speaking, of global fields) provide typical examples of integral domains with finite quotients.

2. As $\fD [R]$ is commutative, it can be identified with the C*-algebra of continuous complex-valued functions on its spectrum. But it turns out that $\Spec (\fD [R])$ is the profinite completion of $R$. Thus, for rings of integers, we get the maximal compact subrings of the finite adeles, which are important objects in number theory.

So much for the finite case, let us turn to general rings now.

\subsection{Towards the general case}

As we have pointed out, the finiteness condition 
\bgloz
  \card{ R / (b) } < \infty \fa b \in R, b \neq 0
\egloz
was heavily used in the formulation of Relation II. If we look at more general rings, this relation will not make sense any more as it is impossible to sum up infinitely many projections.

But at the same time, Relation II contains precious information: It more or less completely determines the structure of the C*-subalgebra $ \fD [R] $ (see \cite{CuLi}, 3.1). Since we do not want to lose all the information, we have to look for an alternative way of describing this commutative C*-subalgebra. 

To this end, let us have a look at the operators $U^a, S_b$ on $ \ell^2 (R) $ again:

As we have stated above, $U^a S_b S_b^* U^{-a} \xi_r = \chf{a+(b)} (r) \xi_r$. 

Let us denote the operator $\xi_r \loma \chf{X} (r) \xi_r$ by $\E{X}$, where $X$ is any subset of $R$. Then, we immediately deduce
\bglnoz
  && \E{a+(b)} \cdot \E{c+(d)} = \E{(a+(b)) \cap (c+(d))}, \\
  && \E{\underset{a+(b) \in R/(b)}{\dotcup} a+(b)} = \E{R} = \sum_{a+(b) \in R/(b)} \E{a+(b)}.
\eglnoz
The crucial observation is that we are given a projection-valued, finitely addivitive spectral measure on the smallest (set-theoretical) algebra over $R$ generated by cosets of all prinicipal ideals. The idea is to view this spectral measure as additional data, on the one hand independent from the unitaries and isometries, but on the other hand allowing an interaction with the original generators and thereby incorporating the ring structure of $R$.

\section{Ring C*-algebras}
\label{R}
Now, we associate to an arbitrary ring with unit two C*-algebras: the reduced and the full ring C*-algebra. However, our construction requires as additional input the choice of a set-theoretical algebra over our ring. Then, we point out how to choose this algebra compatible with the ideal strucure of our ring. Finally, we will see that for such a compatible choice, our construction really generalizes the one of \cite{CuLi}.

Let $R$ be an arbitrary unital ring and denote by $R \reg$ the set of regular elements of $R$: $R \reg = \menge{r \in R}{r \text{ is not a zero-divisor}}$. Moreover, take any set-theoretical algebra $\cC$ over $R$, by which we mean a family of subsets of $R$ containing $R$ and closed under finite unions, finite intersections and complements. Additionally, we require $\cC$ to be invariant under injective affine transformations, which means: $X \in \cC, a \in R, b \in R \reg \Rarr a+b \cdot X \in \cC$. For brevity, we will say that $\cC$ is $P_R$-invariant if $\cC$ satisfies this last property.

Again, the reduced C*-algebra of $R$ is given by concrete operators on $ \ell^2 (R) $:

\bdefin
The families of operators on $ \ell^2 (R) $ given by
\bglnoz
  \E{X} \xi_r &=& \chf{X} (r) \xi_r \text{ for } X \in \cC \\
  U^a \xi_r &=& \xi_{a+r} \text{ for } a \in R \text{ and} \\
  S_b \xi_r &=& \xi_{b \cdot r} \text{ for } b \in R \reg
\eglnoz
give rise to the reduced ring C*-algebra 
\bgloz
  \fA_r [R,\cC] \defeq C^* \rukl{ \menge{ \E{X}, U^a, S_b }{ X \in \cC, a \in R, b \in R \reg } }.
\egloz
\edefin

The full ring C*-algebra is again given as a universal C*-algebra:

\bdefin
We define the full ring C*-algebra of $R$ with respect to $\cC$, denoted by $\fA [R,\cC]$, as the universal C*-algebra generated by
\bglnoz
  && \text{projections} \menge{\e{X}}{X \in \cC}, \\
  && \text{unitaries } \menge{u^a}{a \in R} \\
  && \text{and isometries } \menge{s_b}{b \in R \reg}
\eglnoz 
satisfying

I. The action of $P_R$ on $\cC$ via affine transformations is implemented by the semigrouphomomorphism $ P_R \ni (a,b) \loma u^a s_b $

II. The map $ \cC \ni X \loma \e{X} $ defines a finitely additive spectral measure.
\edefin

\bremark
Let us write out what Relations I. and II. mean. The generators should satisfy
\bglnoz
  && \text{I.(i) } u^a s_b u^c s_d = u^{a+bc} s_{bd} \\
  && \text{I.(ii) } \Ad(u^a s_b) (\e{X}) = \e{a+b \cdot X}
\eglnoz
and 
\bglnoz
  && \text{II.(i) } \e{R} = 1 \\
  && \text{II.(ii) } \e{X} \cdot \e{Y} = \e{X \cap Y} \\
  && \text{II.(iii) } \e{X \dotcup Y} = \e{X} + \e{Y}.
\eglnoz
\eremark

One immediately sees that $\fA [R,\cC]$ exists because all the generators have norm less or equal to $1$. And by universal property of $\fA [R,\cC]$, we always have the canonical projection $\pi: \fA [R,\cC] \lori \fA_r [R,\cC]$ sending generators to generators as it follows directly from a computation that the generators of $\fA_r [R,\cC]$ satisfy I. and II. as well. In particular, $\fA [R,\cC]$ is not trivial.

We only consider regular elements because of two reasons: First of all, the operator $S_b$ would not define a bounded operator in general if $b$ was a zero-divisor, and it is not clear how to modify $S_b$ (for instance by choosing a smaller support projection) to solve this problem. Secondly, restricting to regular elements gives a nice description of the ring C*-algebras associated to direct products (see Proposition \ref{prod}). This result would be destroyed if one wanted to consider zero-divisors as well.

From our motivation, this looks like a reasonable generalization of the construction which we had in the finite case. And indeed, we will show that this new construction really extends the former one in a satisfactory way. But first of all, let us clarify the role of $\cC$.

\subsection{Natural choices of $\cC$}
Since the idea is that our ring C*-algebras should incorporate the ring structure of $R$, it is natural to start with an arbitrary family $\cF$ of right ideals of $R$ and to consider the smallest $P_R$-invariant algebra $\cC(\cF)$ over $R$ generated by $\cF$. Let us denote the full ring C*-algebra $\fA [R,\cC(\cF)]$ by $\fA_{\cF} [R]$ and the reduced version $\fA_r [R,\cC(\cF)] $ by $\fA_{r,\cF} [R]$.

In this situation, typical subsets of $R$ in $ \cC(\cF) $ are given by cosets of the form 
\bgloz
  a + \bigcap_{i=1}^n (b_i \cdot I_i) \text{ with } b_i \in R \reg, I_i \in \cF \cup \gekl{R}.
\egloz
This can be stated more precisely as follows:

\blemma
\label{span}
Let 
\bgloz
  \cF' = \menge{\bigcap_{i=1}^n (b_i \cdot I_i)}{n \in \Zz \pos, b_i \in R \reg, I_i \in \cF \cup \gekl{R}}.
\egloz
Then, $\lspan \rukl{ \menge{\e{X}}{X \in \cC(\cF)} } = \lspan \rukl{ \menge{\e{a+I}}{a \in R, I \in \cF'} }$.
\elemma

For the proof, we need the following observation

\blemma
\label{int}
Let $G_i$ be subgroups of a group $(G,+)$ and take elements $g_i \in G$ ($i=1,2$). Then, we either have that the intersection $(g_1 + G_1) \cap (g_2 + G_2)$ is empty or of the form $g + (G_1 \cap G_2)$ for some $g \in G$.
\elemma

\bproof
If $(g_1 + G_1) \cap (g_2 + G_2) \neq \emptyset$, then take any $g \in (g_1 + G_1) \cap (g_2 + G_2)$. Then, $g + G_1 = g_1 + G_1$ and $g + G_2 = g_2 + G_2$ and thus, 
\bgloz
  (g_1 + G_1) \cap (g_2 + G_2) = (g + G_1) \cap (g + G_2) = g + (G_1 \cap G_2).
\egloz
\eproof

\bproof[Proof of Lemma \ref{span}]
Let us denote $ \lspan \rukl{ \menge{\e{a+I}}{a \in R, I \in \cF'} } $ by $ \fD $. 
To show the inclusion \an{$ \subseteq $} (the other one is obvious), first of all note that $ \e{I} \in \fD $ for all $ I \in \cF $. Moreover, the family 
$\cC' \defeq \menge{X \subseteq R}{\e{X} \in \fD}$ is an $P_R$-invariant algebra over $R$ as $\cC'$ is closed under \dotso

\dotso finite intersections: 
By Relation II.(ii), it suffices to show that $\fD$ is multiplicatively closed. Now, for any $a, \ti{a} \in R; I, \ti{I} \in \cF'$, 
the intersection $(a+I) \cap (\ti{a}+\ti{I})$ is of the form $r + (I \cap \ti{I})$ by Lemma \ref{int}. 
As $\cF'$ is closed under finite intersections by construction, we have $\e{a+I} \cdot \e{\ti{a}+\ti{I}} = \e{r+(I \cap \ti{I})} \in \fD$. 
Thus, $\fD$ is multiplicatively closed.

\dotso complements: This follows from $\e{X \co} = 1 - \e{X}$ (see II.(i) and II.(iii)).

\dotso finite unions: This is a consequence of the first two properties.

\dotso injective affine transformations: This is clear by definition of $\fD$ and $\cF'$.

Hence it follows that $\cC(\cF) \subseteq \cC'$, which implies $\menge{\e{X}}{X \in \cC(\cF)} \subseteq \fD$ and thus, we have shown \an{$\subseteq$}.
\eproof

\bcor
\label{gen}
$\fA_\cF [R]$ is generated by $\menge{\e{I}, u^a, s_b}{I \in \cF, a \in R, b \in R \reg}$.
\ecor

\bproof
This is an immediate consequence of Lemma \ref{span} as for any $X \in \cC(\cF)$ we have $\e{X} \in \lspan(\menge{\e{a+I}}{a \in R, I \in \cF'})=\lspan(\menge{u^a \e{I} u^{-a}}{a \in R, I \in \cF'})$, and for $I = \bigcap_{i=1}^n (b_i \cdot I_i) \in \cF' \text{ with } b_i \in R \reg, I_i \in \cF \cup \gekl{R}$, we have:
\bgloz
  \e{I} = \prod_{i=1}^n \e{b_i \cdot I_i} = \prod_{i=1}^n s_{b_i} \e{I_i} s_{b_i}^* \in C^* \rukl{ \menge{ \e{I}, u^a, s_b }{ I \in \cF, a \in R, b \in R \reg } }.
\egloz
Thus, for any $ X \in \cC ( \cF ) $, $ \e{X} $ lies in $C^*(\menge{\e{I}, u^a, s_b}{I \in \cF, a \in R, b \in R \reg})$.
\eproof

Since we always have the surjection $\pi: \fA [R,\cC] \lori \fA_r [R,\cC]$, all these results hold for the reduced C*-algebras as well. 

\bremark
\label{fAR}

In case we choose $ \cF = \emptyset $, Corollary \ref{gen} tells us that 
\bgloz
  \fA_{\emptyset} [R] = C^* \rukl{ \menge{ u^a, s_b }{ a \in R, b \in R \reg } }.
\egloz
Therefore, in this case $ \fA_{\emptyset} [R] $ is completely determined by the ring itself, without extra choices of ideals or subsets. This then justifies to drop the subscript and to denote this particular ring C*-algebra by $ \fA [R] $. We do so in the reduced case as well ($ \fA_r [R] \defeq \fA_{r,\emptyset} [R] $).
\eremark

\subsection{Compatibility}
To see that the construction we just introduced really generalizes our former one, we observe that both constructions coincide in the case of integral domains with finite quotients, no matter which family of ideals one chooses (as long as one excludes the trivial ideal).

\blemma
\label{com}

Let $R$ be an integral domain with finite quotients. For any family $\cF$ of nontrivial ideals of $R$, the natural map $\fA [R] \lori \fA_{\cF} [R]$ sending generators to generators exists and is an isomorphism.
\elemma

\bproof
This map exists because of the universal property of $ \fA [R] $ (see Definition \ref{fin}) and since in $ \fA_{ \cF } [R] $, the equation
\bgloz
  1 = \sum_{a+(b) \in R/(b)} u^a \e{b} u^{-a} 
\egloz
holds as well because of Relation II.(iii) and the fact that 
\bgloz
  R = \dotcup_{a+(b) \in R/(b)} (a+(b)).
\egloz
Let us first assume that $R$ is not a field. In this case, we know that $ \fA [R] $ is simple (see \cite{CuLi}, 3.3, Theorem 1). Therefore, the map above is injective, as it is certainly nonzero. It remains to prove surjectivity.

Now, the point is that no matter which family $ \cF $ of ideals we choose, we will always end up with
\bgloz
  \cC ( \cF ) = \menge{ \bigcup_{i=1}^n (a_i + I_i)}{ I_i \text{ nontrivial ideal of } R, a_i \in R }.
\egloz
This can be seen as follows:

Denote the right hand side by $ \ti{ \cC } $. We know $ \cF \subseteq \ti{ \cC } $ and $ R \in \ti{ \cC } $. Moreover, $ \ti{ \cC } $ is closed under \dotso

\dotso finite intersections because of Lemma \ref{int}.

\dotso complements since 
\bgloz
  (\bigcup_i (a_i + I_i)) \co = \bigcap_i (a_i + I_i) \co \text{ and } (a_i + I_i) \co = \bigcup_{a_j + I_i \neq a_i + I_i} (a_j + I_i)
\egloz
where the union is finite because $R$ has finite quotients by assumption.

\dotso finite unions by definition.

\dotso injective affine transformations by definition.

Hence, we deduce $ \cC ( \cF ) \subseteq \ti{ \cC } $. 

To see the other inclusion, take any nontrivial ideal $I$ of $R$ and some $ 0 \neq b \in I $. By assumption, we have $ \card{ I/(b) } < \infty $. This implies that $I$ can be written as a finite union
\bgloz
  I = \bigcup_{c+(b) \in I/(b)} (c+(b)). 
\egloz
Now, we can deduce that $I$ lies in $ \cC ( \cF ) $ since $ (b) = b \cdot R $ lies in $ \cC (\cF) $ and because $ \cC (\cF) $ is closed under additive translations and finite unions. Thus, $ \ti{\cC} \subseteq \cC(\cF) $. So we have seen our claim that for any family $ \cF $ of nontrivial ideals of $R$, $\cC ( \cF )$ and $ \ti{ \cC } $ coincide. 

Therefore, we can take $ \cF = \emptyset $, and then we know that $ \fA_{\cF} [R] $ is generated by $\menge{ u^a, s_b }{ a \in R, b \in R \reg }$ by Remark \ref{fAR}. This shows surjectivity and thus our claim.

Finally, if $R$ is a field, our constructions will both yield the maximal group C*-algebra of the $ax+b$-group $P_R$. Therefore, our result holds here as well.
\eproof

Actually, it is possible to deduce Lemma \ref{com} without using simplicity of $ \fA [R] $. Instead, we can write both C*-algebras as crossed products of commutative C*-algebras by the $ax+b$-semigroup $P_R$ (see Proposition \ref{cropro}). Then, it remains to identify the spectra of these commutative C*-subalgebras in a $P_R$-equivariant way (compare Remark \ref{equi}), and this can be done by describing these spectra as inverse limits (in the spirit of \cite{CuLi}, 4.2, Observation 1).

\section{Basic properties}
\label{B}

Let us derive some immediate consequences from the definitions. It turns out that ring C*-algebras can be described as semigroup crossed products. Furthermore, the C*-algebra associated to a product of rings can be identified with the tensor product of the C*-algebras associated to each of the factors.

\subsection{Crossed product description}

As in the finite case, $\fA [R,\cC]$ admits a description as a crossed product. But since $R \reg$ is only a semigroup in the general case, we have to consider crossed products by semigroups. The basics of this theory are explained in the appendix.

First of all, there is a canonical commutative C*-subalgebra of $ \fA [R,\cC] $, namely
\bgloz
  \fD [R,\cC] \defeq C^* \rukl{ \menge{\e{X}}{X \in \cC} } \subseteq \fA [R,\cC].
\egloz
Commutativity of $ \fD [R,\cC] $ follows from Relation II.(ii). 

Moreover, Relation I yields an action of $P_R$ on $ \fD [R,\cC] $ given by the semigrouphomomorphism $\alpha: P_R \lori \End (\fD [R,\cC]); (a,b) \loma \Ad (u^a s_b)$. 

Now, the following observation merely reformulates the definition of $\fA [R,\cC]$:

\bprop
\label{cropro}

$\fA [R,\cC]$ is isomorphic to $\fD [R,\cC] \rte_{\alpha} P_R$.
\eprop

\bproof
Just compare the universal properties of these C*-algebras. Both of them are defined as universal C*-algebras generated by $\fD [R,\cC]$ and isometries 
$\menge{V_{(a,b)}}{(a,b) \in P_R}$ ($u^a s_b$ for $\fA [R,\cC]$) with $\Ad (V_{(a,b)}) \vert_{\fD [R,\cC]} = \alpha (a,b)$ ($\Ad(u^a s_b) = \alpha (a,b)$ holds in $\fA [R,\cC]$ by definition).
\eproof

As we will see, this simple observation already has some consequences on the structure of these ring C*-algebras, at least if we consider commutative rings.

\subsection{Direct products of rings}

In analogy to the case of groups, we have the following

\bprop
\label{prod}
Let $R_i$ be two unital rings together with $P_{R_i}$-invariant algebras $\cC_i$ over $R_i$ ($i=1,2$). Then, 
\bglnoz
  && \fA_r [R_1 \times R_2, \cC(\cC_1 \times \cC_2)] \cong \fA_r [R_1, \cC_1] \otimes_{\min} \fA_r [R_2, \cC_2] \text{ and} \\
  && \fA [R_1 \times R_2, \cC(\cC_1 \times \cC_2)] \cong \fA [R_1, \cC_1] \otimes_{\max} \fA [R_2, \cC_2],
\eglnoz
where $ \cC(\cC_1 \times \cC_2) $ is the smallest $ P_{R_1 \times R_2} $-invariant algebra over $R_1 \times R_2$ containing $\cC_1 \times \cC_2 = \menge{X_1 \times X_2 \subseteq R_1 \times R_2}{X_i \in \cC_i}$.
\eprop

\bproof
Let us write $\cC$ for $\cC(\cC_1 \times \cC_2)$. First of all, it is helpful to observe that $(R_1 \times R_2) \reg = R_1 \reg \times R_2 \reg$. This holds true since multiplication on $R_1 \times R_2$ is defined componentwise. 

Now, in the reduced case, $ \fA_r [R_1,\cC_1] \otimes_{\min} \fA_r [R_2,\cC_2] $ is the closure of the algebraic tensor product 
$ \fA_r [R_1,\cC_1] \odot \fA_r [R_2,\cC_2] $ in $ \cL (\ell^2(R_1) \otimes \ell^2(R_2)) $. This is just the definition of the minimal tensor product. Moreover, consider the unitary
\bgloz
  W: \ell^2(R_1 \times R_2) \lori \ell^2(R_1) \otimes \ell^2(R_2); \xi_{(r_1,r_2)} \loma \xi_{r_1} \otimes \xi_{r_2}.
\egloz
It satisfies 
\bgln
\label{propW}
  W \E{X_1 \times X_2} W^* &=& \E{X_1} \otimes \E{X_2} \\
  W U^{(a_1,a_2)} W^* &=& U^{a_1} \otimes U^{a_2} \nonumber \\
  W S_{(b_1,b_2)} W^* &=& S_{b_1} \otimes S_{b_2} \nonumber
\egln
for all $X_i \in \cC_i$, $a_i \in R_i$ and $b_i \in R_i \reg$.

Thus, $W (\fA_r [R_1 \times R_2, \cC(\cC_1 \times \cC_2)]) W^* \subseteq \fA_r [R_1, \cC_1] \otimes_{\min} \fA_r [R_2, \cC_2]$. The other inclusion follows from
\bglnoz
  && \fA_r [R_1, \cC_1] \otimes_{\min} \fA_r [R_2, \cC_2] = C^* \rukl{ \fA_r [R_1, \cC_1] \otimes 1, 1 \otimes \fA_r [R_2, \cC_2] } \\
  &\text{and}& \fA_r [R_1, \cC_1] \otimes 1, 1 \otimes \fA_r [R_2, \cC_2] \subseteq W (\fA_r [R_1 \times R_2, \cC(\cC_1 \times \cC_2)]) W^* 
\eglnoz
where we used the properties of the unitary $W$ listed under (\ref{propW}).

To identify the full ring C*-algebras, we compare the universal properties:

$\fA [R_1 \times R_2, \cC(\cC_1 \times \cC_2)]$ is generated by
\bglnoz
  && \menge{\e{X_1 \times R_2}, u^{(a_1,0)}, s_{(b_1,1)}}{X_1 \in \cC_1, a_1 \in R_1, b_1 \in R_1 \reg} \\
  &\text{and}& \menge{\e{R_1 \times X_2}, u^{(0,a_2)}, s_{(1,b_2)}}{X_2 \in \cC_2, a_2 \in R_2, b_2 \in R_2 \reg}, 
\eglnoz
and the relations are just the defining relations for $ \fA [R_i,\cC_i] $ for each of these families of generators. Moreover, these two families of generators commute with each other. For instance, 
\bgloz
  s_{(b_1,1)} \e{R_1 \times X_2} s_{(b_1,1)}^* 
  = \e{(b_1 \cdot R_1) \times X_2}
  = \e{R_1 \times X_2} \cdot \e{(b_1 \cdot R_1) \times R_2}
  = \e{R_1 \times X_2} \cdot s_{(b_1,1)} s_{(b_1,1)}^*
\egloz
which implies $\eckl{s_{(b_1,1)}, \e{R_1 \times X_2}}=0$. But $\fA [R_1, \cC_1] \otimes_{\max} \fA [R_2, \cC_2]$ admits exactly the same description. Hence these C*-algebras are isomorphic.
\eproof

\section{Commutative Rings}
\label{C}

In the case of commutative rings, we can deduce rather strong results on the inner structure of the associated ring C*-algebras. It turns out that the C*-algebras of commutative rings always lie in the nuclear UCT class and that one can give necessary and sufficient conditions for them to be purely infinite and simple (see \cite{Ror} for the definition of \an{purely infinite and simple}). 

\subsection{Nuclearity and UCT}

Let $R$ be a commutative ring with unit and let $\cC$ be a $P_R$-invariant algebra over $R$. 

First of all, $\fA [R,\cC]$ is isomorphic to $\fD [R,\cC] \rte_{\alpha} P_R$ as shown in Proposition \ref{cropro}. This can be used to derive some properties of our constructions.

As pointed out in \cite{La}, crossed products by semigroups and ordinary crossed products are closely related (at least in nice cases). In our case, the $ax+b$-semigroup $P_R$ is an Ore semigroup acting on $ \fD [R,\cC] $ via injective endomorphisms, which is a particularly nice situation.

We obtain the associated ordinary C*-dynamical system by formally inverting $ \alpha (a,b) $ for all $ (a,b) \in P_R $: Consider the inductive system given by 
\begin{itemize}
  \item the C*-algebras $ \fD (R,\cC)_{ (a,b) } = \fD [R,\cC] $ for all $ (a,b) \in P_R $
  \item the structure maps $ \alpha(c,d) : \fD (R,\cC)_{ (a,b) } \lori \fD (R,\cC)_{ (a,b) \cdot (c,d) } $.
\end{itemize}
The inductive limit $ \fD (R,\cC) $ carries a natural action $ \overline{ \alpha } $ of $ P_R^{-1} \cdot P_R = P_{ Q(R) } $ by automorphisms. Here, $Q(R)$ is the ring of fractions $(R \reg)^{-1} R$. $ \overline{ \alpha } $ extends $ \alpha $ with respect to the natural embedding $ i: \fD [R,\cC] \into \fD (R,\cC) $ ($i$ is given by the composition $ \fD [R,\cC] \lori \fD (R,\cC)_{ (0,1) } \lori \fD (R,\cC) $). In Laca's notation, the C*-dynamical system $(\fD (R,\cC) , P_{ Q(R) } , \overline{ \alpha })$ is the minimal automorphic dilation of $(\fD [R,\cC] , P_R , \alpha)$ (see \cite{La}, Theorem 2.1.1.).

It turns out that the corresponding crossed product algebras are strongly Morita equivalent:

\blemma
\label{cor}

There exists an isomorphism 
\bgloz
  \Phi : \fA [R,\cC] \cong \fD [R,\cC] \rte_{\alpha} P_R \lori i(1) (\fD (R,\cC) \rtimes_{\overline{\alpha}} P_{Q(R)}) i(1) 
\egloz
satisfying $\Phi \vert_{\fD [R, \cC]}=i$. 

Moreover, $\Phi \lru u^a s_b \rru = i(1) U_{(a,b)} i(1)$ where $U: P_{Q(R)} \ri \cU \cM(\fD (R,\cC) \rtimes_{\overline{\alpha}} P_{Q(R)})$ is the associated unitary representation of $P_{Q(R)}$. 

Furthermore, $ i(1) $ is full projection in $ \fD (R,\cC) \rtimes_{ \overline{ \alpha } } P_{ Q(R) } $.
\elemma

\bproof
Apply \cite{La}, Theorem 2.2.1.
\eproof

Let us derive two consequences from this observation:

- $ \fA [R,\cC] $ is in the UCT class:

As $ \fD [R,\cC] $ is Morita equivalent to $ \fD (R,\cC) \rtimes_{ \overline{ \alpha } } P_{ Q(R) } $, it suffices to consider the latter C*-algebra. 
But $ \fD (R,\cC) $ is a commutative C*-algebra, so that $ \fD (R,\cC) \rtimes_{ \overline{ \alpha } } P_{ Q(R) } $ can be written as the corresponding groupoid C*-algebra. The transformation groupoid is amenable since $ P_{ Q(R) } $ is solvable, and hence amenable. Therefore, a general result of \cite{Tu} shows that 
$ \fD (R,\cC) \rtimes_{ \overline{ \alpha } } P_{ Q(R) } $ satisfies UCT.

- $ \fA [R,\cC] $ is nuclear:

As $ P_{ Q(R) } $ is amenable and $ \fD (R,\cC) $ is abelian, thus nuclear, $ \fD (R,\cC) \rtimes_{ \overline{ \alpha } } P_{ Q(R) } $ is nuclear (see \cite{Ror}). Since hereditary C*-subalgebras of nuclear C*-algebras are nuclear again (see \cite{Ror}), we conclude that
\bgloz
  \fA [R,\cC] \cong i(1) (\fD (R,\cC) \rtimes_{ \overline{ \alpha } } P_{ Q(R) })i(1) 
\egloz
must be nuclear again.

\subsection{Purely infinite and simple C*-algebras}

Let $\cF$ be a family of ideals in $R$. Moreover, let $Z$ denote the set of zero-divisors in $R$.

Our present goal is to investigate under which conditions $\fA_{\cF} [R]$ is purely infinite and simple. First of all, we have to impose two technical conditions on our ring, namely:
\bglnoz
  && \bigcap_{b \in R \reg} (b) = (0), \\
  && \text{for all ideals } I \subseteq Z \text{ in } R \text{ and for all } b \in R \reg, \card{ \qr{(b)}{(b) \cap I} } = \infty.
\eglnoz
Our main result is the following:

\btheooz
Assume that $R$ is a commutative ring with unit satisfying the two conditions above. Then, for any family $\cF$ of ideals in $R$, the ring C*-algebra $ \fA_{\cF} [R] $ is purely infinite and simple if and only if all $I \in \cF$ satisfy $I \nsubseteq Z$.
\etheooz

\subsubsection{A criterion}

The central idea is to describe a particular faithful conditional expectation by sufficiently small projections. This idea already appears in \cite{Cun2}. Since then, it has been continuously modified to detect purely infinite and simple C*-algebras in various other situations (see, for example, \cite{ExLa} or \cite{LaSp}). But the basic idea has remained the same. So, we would like to start with this abstract concept. 

\bprop
\label{cri}

Let $\cA$ be a dense *-algebra of a unital C*-algebra $A$. Assume that $\Theta$ is a faithful conditional expectation on $A$ such that for every $0 \neq x \in \cA_+$ there exist finitely many projections $f_i \in A$ with
\bgln
\label{c1}
  && f_i \perp f_j \fa i \neq j; f_i \overset{s_i}{\sim} 1 \text{ via isometries } s_i \in A \fa i, \\
\label{c2}
  && \norm{\sum_i f_i \Theta(x) f_i } = \norm{\Theta(x)}; f_i x f_i = f_i \Theta(x) f_i \in \Cz f_i \fa i.  
\egln
Then $A$ is purely infinite and simple.
\eprop

\bproof
Take $0 \neq a \in A$. Without loss of generality we can assume $a \in A_+$ and $\norm{\Theta (a)} = 1$. As $\cA$ is a dense *-subalgebra, we can find $x \in \cA_+$ with the properties $\norm{a-x} < \halb$ and $\norm{\Theta(x)} = 1$.

By hypothesis, we can find projections $f_i$ for this element $x$ with 
\bgloz
  1 = \norm{ \Theta (x) } = \norm{ \sum_i f_i \Theta (x) f_i } \text{ (see (\ref{c2}))}.
\egloz

As $ x $ is positive, $ f_i \Theta (x) f_i $ is positive, too. Thus, we can write $f_i \Theta (x) f_i = \lambda_i f_i$ with $\lambda_i \geq 0$ by (\ref{c2}). Now, using mutual orthogonality of the projections (compare (\ref{c1})), we can conclude that 
\bgloz
  \sup_i \lambda_i = \norm{ \sum_i \lambda_i f_i } = \norm{ \sum_i f_i \Theta (x) f_i } = 1.
\egloz 
Hence it follows that among these projections $f_i$, there is a particular one, say $f$, with $ f = f \Theta (x) f = f x f $ (which means $\lambda=1$ for this projection). By (\ref{c1}), there exists an isometry $ s \in A $ with $ s s^* = f $.

Therefore, we can calculate
\bglnoz
  && \norm{ s^* a s - 1 } = \norm{ s^* a s - s^* f s } = \norm{ s^* a s - s^* f x f s } \\
  &=& \norm{ s^* a s - s^* x s } \leq \norm{ a-x } < \halb.
\eglnoz
This implies that $ s^* a s $ is invertible in $ A $. Thus, setting $ y \defeq (s^* a s)^{-1} s^* $, $ z = s $ we have $ y a z = 1 $.
\eproof

\subsubsection{Preparations}

Let us fix a commutative ring with unit, say $R$, for the rest of this section. Moreover, let $\cF$ be some family of ideals in $R$. 

To apply Proposition \ref{cri} to our situation ($A = \fA_{\cF} [R]$), we have to construct a faithful conditional expectation together with a suitable dense *-subalgebra.

\blemma
\label{dense}
Let 
\bgloz
  \cF' = \menge{ \bigcap_{i=1}^n (b_i \cdot I_i)}{b_i \in R \reg, I_i \in \cF \cup \gekl{R}}
\egloz
as in Lemma \ref{span}.

The *-subalgebra of $\fA_{\cF} [R]$ generated by $\menge{\e{I},u^a,s_b}{I \in \cF, a \in R, b \in R \reg}$ coincides with $\cA_{\cF} [R] \defeq \lspan(\menge{s_b^* u^{-a} \e{I} u^{a'} s_{b'}}{I \in \cF', a \in R, b \in R \reg})$.
\elemma

\bproof
All we have to show is that $\cA_{\cF} [R]$ is multiplicatively closed since 
\bglnoz
  && \menge{\e{I},u^a,s_b}{I \in \cF, a \in R, b \in R \reg} \\
  &\subseteq& \cA_{\cF} [R] \\
  &\subseteq& \text{ *-} \alg \rukl{ \menge{\e{I},u^a,s_b}{I \in \cF, a \in R, b \in R \reg} }.
\eglnoz
The following computation shows that $ \cA_{\cF} [R] $ is multiplicatively closed:
\bglnoz
  && \lru s_b^* u^a \e{I} u^{-a'} s_{b'} \rru \lru s_d^* u^c \e{J} u^{-c'} s_{d'} \rru \\
  &=& s_b^* u^{a-a'} \underbrace{ \lru u^{a'} \e{I} u^{-a'} \rru }_{ = \e{a'+I} } s_{b'} s_{b'}^* 
  s_d^* s_{b'} \underbrace{ \lru u^c \e{J} u^{-c} \rru }_{ = \e{c+J} } u^{c-c'} s_{d'} \\
  &=& s_b^* u^{a-a'} s_d^* \underbrace{ \lru s_d \e{a'+I} s_d^* \rru }_{ = \e{d(a'+I)}  } \underbrace{ \lru s_d \e{(b')} s_d^* \rru }_{ = \e{(db')} } 
  \underbrace{ \lru s_{b'} \e{c+J} s_{b'}^* \rru }_{ = \e{b'(c+J)} } s_{b'} u^{c-c'} s_{d'} \\
  &=& s_b^* u^{a-a'} s_d^* \e{ d(a'+I) \cap (db') \cap b'(c+J) } s_{b'} u^{c-c'} s_{d'} \\
  &=& s_{bd}^* u^{d(a-a')} \e{ d(a'+I) \cap (db') \cap b'(c+J) } u^{b'(c-c')} s_{b'd'}
\eglnoz
for any $ I, J \in \cF' $, $ a, a', c, c' \in R $ and $ b, b', d, d' \in R \reg $. 

Now, by Lemma \ref{int}, the intersection $ d(a'+I) \cap (db') \cap b'(c+J) $ is either empty or of the form $ \ti{a} + (d \cdot I) \cap (db') \cap (b' \cdot J) $. If it is empty, the product above will vanish, hence it will lie in $ \cA_{\cF} [R] $. If the intersection is not empty, then the product will be 
\bglnoz
  && s_{bd}^* u^{d(a-a')} \e{ \ti{a} + (d \cdot I) \cap (db') \cap (b' \cdot J) } u^{b'(c-c')} s_{b'd'} \\
  &=& s_{bd}^* u^{d(a-a')+\ti{a}} \e{ (d \cdot I) \cap (db') \cap (b' \cdot J) } u^{-\ti{a}+b'(c-c')} s_{b'd'}
\eglnoz
which will again lie in $ \cA_{\cF} [R] $.
\eproof

\blemma
\label{fce}
With the same notations as in Lemma \ref{cor}, there exists a faithful conditional expectation $\Theta: \fA_{\cF}[R] \lori \Phi^{-1} (i(1)(\fD_{\cF}(R))i(1)) \subseteq \fA_{\cF}[R]$ with 
\bgl
  \label{condexp}
  \Theta (s_b^* u^{-a} \e{I} u^{a'} s_{b'}) = \delta_{a,a'} \delta_{b,b'} s_b^* u^{-a} \e{I} u^{a} s_{b}
\egl
for all $ I \in \cF' $, $ a, a' \in R $ and $ b, b' \in R \reg $.
\elemma

\bproof
We know that there is a faithful conditional expectation 
\bgloz
  \overline{ \Theta } : \fD_{ \cF } (R) \rtimes_{ \overline{ \alpha } } P_{ Q(R) } \lori \fD_{ \cF } (R) 
\egloz
which satisfies
\bgloz
  \overline{\Theta} (U_{(a,b)}^* x U_{(a',b')}) = \delta_{a,a'} \delta_{b,b'} U_{(a,b)}^* x U_{(a,b)}
\egloz
for all $x \in \fD_{\cF} (R)$ (for instance, as $P_{Q(R)} = Q(R) \rtimes Q(R) \reg$, we can apply the construction described in \cite{Bla}, II.10.4.17 iteratively to the dual groups of $Q(R) \reg$ and $Q(R)$).

It immediately follows that the composition $ \Theta \defeq \Phi^{-1} \circ \overline{ \Theta } \circ \Phi $ is a faithful conditional expectation satisfying (\ref{condexp}).
\eproof

\subsubsection{A general fact on subgroups}

The following fact will be useful:

\blemma
Assume that $G_i$ are subgroups of an abelian group $(G,+)$ with $ \card{G/G_i} = \infty $ for $1 \leq i \leq n$. Then, for all $g_i \in G$, we have 
\bgloz
  G \neq \bigcup_{i=1}^n (g_i+G_i).
\egloz
\elemma

\bproof
Let $m=\# \gekl{G_1, \dotsc, G_n} $ be the number of pairwise distinct subgroups among the $G_i$. Without loss of generality, we can assume that the subgroups are indexed so that $\gekl{G_1, \dotsc, G_n} = \gekl{G_1, \dotsc, G_m}$. We prove our assertion inductively.

$m=1$: The claim follows from $ \card{G/G_i} = \infty $.

$m>1$: Assume that we have proven the claim for any $m-1$ pairwise distinct subgroups, and assume that 
\bgloz
  G = \bigcup_{i=1}^n (g_i+G_i).
\egloz
There are two possible cases:

1. There exists $1<j \leq m$ with $\card{\qr{G_1}{G_1 \cap G_j}} < \infty$. Because of $\ql{G_1+G_j}{G_j} \cong \qr{G_1}{G_1 \cap G_j}$, it follows that $\card{ \ql{G_1+G_j}{G_j}} < \infty$. Now, the exact sequence $\ql{G_1+G_j}{G_j} \into G / G_j \onto \qr{G}{G_1+G_j}$, together with $\card{\ql{G_1+G_j}{G_j}} < \infty$, $\card{G / G_j} = \infty$, implies $\card{\qr{G}{G_1+G_j}} = \infty$. 

Define 
$
  \ti{G_i} = 
  \bfa
    G_i \falls G_i \neq G_1, G_i \neq G_j \\
    G_1 + G_j \falls G_i \in \gekl{G_1,G_j} 
  \efa
$. We still have $G = \bigcup_{i=1}^n (g_i+\ti{G}_i)$, but $\# \gekl{\ti{G}_1, \dotsc, \ti{G}_n} \leq m-1$ which contradicts the induction hypothesis.

2. For all $1<j \leq m$, $ \card{ \qr{G_1}{G_1 \cap G_j} } = \infty $. Since we know $ \card{G/G_1} = \infty $, there exists $ g \in G $ with $ g+G_1 \neq g_i+G_i $ for all $ 1 \leq i \leq n $. Then, 
\bglnoz
  g+G_1 &=& (g+G_1) \cap G = \bigcup_{i=1}^n (g+G_1) \cap (g_i+G_i) \\
  &=& \bigcup_{ (g+G_1) \cap (g_i+G_i) \neq \emptyset } (g+G_1) \cap (g_i+G_i) \\
  &=& \bigcup_{ (g+G_1) \cap (g_i+G_i) \neq \emptyset } (\ti{g_i} + (G_1 \cap G_i))
\eglnoz
for some $ \ti{g_i} \in G $, where we used Lemma \ref{int} in the last step. 

This implies
\bgloz
  G_1 = \bigcup_{ (g+G_1) \cap (g_i+G_i) \neq \emptyset } ((\ti{g_i}-g) + (G_1 \cap G_i))
\egloz
and $\ti{g_i}-g$ lies in $G_1$ for all those indices $i$ with $(g+G_1) \cap (g_i+G_i) = \ti{g_i}+(G_1 \cap G_i)$.

But the number of pairwise distinct subgroups of $G_1$ among
\bgloz
  \menge{G_1 \cap G_i}{(g+G_1) \cap (g_i+G_i) \neq \emptyset}
\egloz
has decreased because $ (g+G_1) \cap (g_i+G_i) = \emptyset $ for all $i$ with $G_i = G_1$ by construction. Again, this contradicts the induction hypothesis.
\eproof

\bcor
\label{finun}
Let $I$, $I_i$ be ideals in $R$ with $ \card{ \qr{I}{I \cap I_i} } = \infty $ for all $ 1 \leq i \leq n $. Then, for all $a,a_i \in R$, we have $a+I \nsubseteq \bigcup_{i=1}^n (a_i+(I \cap I_i))$.
\ecor

\bproof
Apply the preceding Lemma to $G=I$, $G_i = I \cap I_i$. Translation by $a \in R$ yields the claim for cosets as well.
\eproof

\subsubsection{Manipulations of projections}
\label{mani}

Our present aim is to show that $ \fA_{\cF} [R] $, together with the dense *-subalgebra $ \cA_{\cF} [R] $ (see Lemma \ref{dense}) and the faithful conditional expectation $\Theta$ (constructed in Lemma \ref{fce}), satisfies our criterion (Proposition \ref{cri}) if $R$ satisfies the two technical conditions we mentioned above and if all $I \in \cF$ satisfy $I \nsubseteq Z$. We still have to find appropriate projections $f_i$ depending on a nontrivial element of $\cA_{\cF} [R]_+$.

To do so, let us consider the following situation: Take $0 \neq x \in \cA_{\cF} [R]_+$. Since $\Theta$ is faithful, we know $\Theta(x) \neq 0$. Moreover, $\Theta(x)$ is a finite sum of the form
\bgloz
  \sum_{(d',X)} \beta_{(d',X)} s_{d'}^* \e{X} s_{d'}
\egloz
where the sum is taken over pairs $(d',X) \in R \reg \times \cC(\cF) $ (by property (\ref{condexp}) of $\Theta$).

Define $d$ as the product of all $d' \in R \reg$ with $\beta_{(d',X)} s_{d'}^* \e{X} s_{d'} \neq 0$ for some $X \in \cC(\cF)$. Then 
$\sum_{(d',X)} \beta_{(d',X)} s_{d'}^* \e{X} s_{d'} = s_d^* (\sum_{(d',X)} \beta_{(d',X)} \e{(d/d') \cdot X}) s_d$. 

Moreover, with $\cF' = \menge{ \bigcap_{i=1}^n (b_i \cdot I_i)}{b_i \in R \reg, I_i \in \cF \cup \gekl{R}}$ as in Lemma \ref{span}, it is possible to write
\bgl
\label{sum}
  \Theta(x) = s_d^* (\sum_{(m,I)} \beta_{(m,I)} \e{m+I})s_d,
\egl
where the sum is taken over pairs $ (m,I) \in R \times \cF' $ with only finitely many nontrivial coefficients (see Lemma \ref{span}).

We can additionally assume $m+I \subseteq (d)$ for all $(m,I)$ with $ \beta_{(m,I)} \neq 0 $, since we can substitute $ \e{m+I} $ by 
$ \e{(d)} \cdot \e{m+I} = \e{(d) \cap (m+I)} $ and $ (d) \cap (m+I) $ is of the form $m'+((d) \cap I)$ if nonempty (see Lemma \ref{int}).

Now, in this situation, we claim the following:

\blemmaoz
Assume that all $I \in \cF$ satisfy $I \nsubseteq Z$. Then, there exist finitely many pairwise orthogonal projections $p_i$ in $\fA_{\cF} [R]$ with
\bgloz
  C^* \rukl{ \gekl{p_i} } = C^* \rukl{ \menge{\e{m+I}}{\beta_{(m,I)} \neq 0} } \text{(see (\ref{sum}))}.
\egloz
Moreover, we can find $b_i' \in R \reg$, $a_i' \in R$ with $ \e{a_i'+(b_i')} \leq p_i \fa i $.
\elemmaoz

The proof of this technical lemma will be broken into several parts:

\blemma
\label{l0}
Let $P_{\cF} [R] = \menge{\e{m'+I'}}{m' \in R, I' \in \cF'}$. We can find finitely many pairwise orthogonal (nontrivial) projections $p_i$ in $ \Zz \text{-} \lspan(P_{\cF} [R]) $ with $C^*(\menge{\e{m+I}}{\beta_{(m,I)} \neq 0}) = C^*(\gekl{p_i})$.
\elemma

\bproof
As the projections commute, we can easily orthogonalize them. And we can show inductively that the coefficients we get by orthogonalizing them are integers.
\eproof

Now, fix a projection $p \in \gekl{p_i}$. As $p$ lies in $ \Zz \text{-} \lspan(P_{\cF} [R]) $, we can write 
\bgl
\label{p1}
  p = \sum_j n_j \e{a_j+I_j} - \sum_{j'} \ti{n}_{j'} \e{\ti{a}_{j'}+\ti{I}_{j'}} 
\egl
with finitely many $n_j$, $\ti{n}_{j'}$ in $\Zz \pos$. 

As a next step, we show

\blemma
\label{l1}
This representation (\ref{p1}) of $p$ can be modified so that for any ideals $J,\ti{J} \in \gekl{I_j, \ti{I}_{j'}}$, the following condition holds: 
\bgl
\label{cond}
  \text{Either } \card{ \qr{J}{J \cap \ti{J}} } = 1 \text{ or } \card{ \qr{J}{J \cap \ti{J}} } = \infty.
\egl
\elemma

\bproof
We will arrange this by the following recursion:

Enumerate the ideals $ \gekl{I_j, \ti{I}_{j'}} $ so that $ \gekl{I_j, \ti{I}_{j'}} = \gekl{J_g} $. Assume that our condition (\ref{cond}) holds for $\gekl{J_1, \dotsc, J_h}$. 
For $h=1$, the condition is automatically satisfied. 

Now, define $J_{h+1}^{(0)} = J_{h+1}$ and for $g=1, \dotsc, h$ 
\bgloz
  J_{h+1}^{(g)} = 
  \bfa
    J_{h+1}^{(g-1)} \falls \card{ \qr{J_{h+1}^{(g-1)}}{J_{h+1}^{(g-1)} \cap J_g} } \in \gekl{1, \infty} \\
    J_g \cap J_{h+1}^{(g-1)} \falls 1 < \card{ \qr{J_{h+1}^{(g-1)}}{J_{h+1}^{(g-1)} \cap J_g} } < \infty
  \efa.
\egloz

We have to argue that this substitution can be done for the corresponding projections in (\ref{p1}) as well: 

If $ \card{ \qr{J_{h+1}^{(g-1)}}{J_{h+1}^{(g-1)} \cap J_g} } = M < \infty $, we can find $r_1, \dotsc, r_M \in R$ with $J_{h+1}^{(g-1)} = (r_1 + (J_{h+1}^{(g-1)} \cap J_g)) \dotcup \dotsb \dotcup (r_M + (J_{h+1}^{(g-1)} \cap J_g))$. Thus, by Relation II.(iii), we can replace $\e{J_{h+1}^{(g-1)}}$ by $\e{r_1 + (J_{h+1}^{(g-1)} \cap J_g)} + \dotsb + \e{r_M + (J_{h+1}^{(g-1)} \cap J_g)}$ in the representation (\ref{p1}) of $p$. This allows us to substitute $J_{h+1}$ by $J'_{h+1} \defeq J_{h+1}^{(h)}$.

We claim: $\card{ \qr{J'_{h+1}}{J'_{h+1} \cap J} } \in \gekl{1, \infty} \fa J \in \gekl{J_1, \dotsc, J_h}$.

Proof of the first claim: To see this, let us prove inductively on $g$ that 
\bgloz
  \card{ \qr{J_{h+1}^{(g)}}{J_{h+1}^{(g)} \cap J} } \in \gekl{1, \infty} \fa J \in \gekl{J_1, \dotsc, J_g}.
\egloz

For $g=1$, the assertion holds by construction.

Assume that we have proven this assertion for $g-1$. Take any $J \in \gekl{J_1, \dotsc, J_g}$. By construction, 
$ \card{ \qr{J_{h+1}^{(g)}}{J_{h+1}^{(g)} \cap J_g} } \in \gekl{1, \infty} $. Thus, it suffices to consider $J \in \gekl{J_1, \dotsc, J_{g-1}}$. 

If $ J_{h+1}^{(g)} = J_{h+1}^{(g-1)} $, we will have $ \card{ \qr{J_{h+1}^{(g)}}{J_{h+1}^{(g)} \cap J} } \in \gekl{1, \infty} $ by induction hypothesis. 

Otherwise, $ J_{h+1}^{(g)} = J_{h+1}^{(g-1)} \cap J_g \subseteq J_{h+1}^{(g-1)} $ and $ \card{ \qr{J_{h+1}^{(g-1)}}{J_{h+1}^{(g-1)} \cap J_g} } < \infty $. 

If $ \card{ \qr{J_{h+1}^{(g-1)}}{J_{h+1}^{(g-1)} \cap J} } = 1 $, we will have $J_{h+1}^{(g)} \subseteq J_{h+1}^{(g-1)} \subseteq J$ which implies 
$ \card{ \qr{J_{h+1}^{(g)}}{J_{h+1}^{(g)} \cap J} } = 1 $.

If $ \card{ \qr{J_{h+1}^{(g-1)}}{J_{h+1}^{(g-1)} \cap J} } = \infty $, consider the exact sequence 
\bgl
\label{ex1}
  \qr{J_{h+1}^{(g)}}{J_{h+1}^{(g)} \cap J} \into \qr{J_{h+1}^{(g-1)}}{J_{h+1}^{(g)} \cap J} \onto J_{h+1}^{(g-1)} / J_{h+1}^{(g)}.
\egl
Now, $\card{ \qr{J_{h+1}^{(g-1)}}{J_{h+1}^{(g)} \cap J} } = \infty$ since $ \card{ \qr{J_{h+1}^{(g-1)}}{J_{h+1}^{(g-1)} \cap J} } = \infty $ and because there is the  projection $\qr{J_{h+1}^{(g-1)}}{J_{h+1}^{(g)} \cap J} \onto \qr{J_{h+1}^{(g-1)}}{J_{h+1}^{(g-1)} \cap J}$. But we also know $\card{ J_{h+1}^{(g-1)} / J_{h+1}^{(g)} } < \infty$. Thus, $ \card{ \qr{J_{h+1}^{(g)}}{J_{h+1}^{(g)} \cap J} } = \infty $ as (\ref{ex1}) is exact. 

Hence, we have seen $\card{ \qr{J_{h+1}^{(g)}}{J_{h+1}^{(g)} \cap J} } \in \gekl{1, \infty}$ for all $J \in \gekl{J_1, \dotsc, J_g}$. This proves our first claim.

As a second step, substitute $J \in \gekl{J_1, \dotsc, J_h}$ by $J' = J \cap J'_{h+1}$ whenever $1 < \card{ \qr{J}{J \cap J'_{h+1}} } < \infty$. Again, Relation II.(iii) allows us to substitute $J$ by $J'$ in the representation (\ref{p1}) of $p$. 

So, we end up with a new set of ideals $\gekl{J'_1, \dotsc, J'_{h+1}}$. We claim: 

For any $J,\ti{J} \in \gekl{J'_1, \dotsc, J'_{h+1}}$, condition (\ref{cond}) holds: $\card{ \qr{J}{J \cap \ti{J} } } \in \gekl{1, \infty}$.

Proof of the second claim:

1. Consider the case $J \in \gekl{J'_1, \dotsc, J'_h}$, $\ti{J} = J'_{h+1}$: By the last step, we have enforced $ \card{ \qr{J}{J \cap \ti{J} } } \in \gekl{1, \infty} $

2. Next, assume $J=J'_{h+1}$, $\ti{J}=J'_g \in \gekl{J'_1, \dotsc, J'_h}$: By construction, we have $J \cap J'_g = J \cap J_g$. Thus, by our first claim, 
\bgloz
  \card{\qr{J}{J \cap \ti{J}}} = \card{\qr{J}{J \cap J_g}} \in \gekl{1, \infty}.
\egloz

3. Finally, let $J=J'_g$ and $\ti{J}=J'_{\ti{g}}$ with $g,\ti{g} \in \gekl{1, \dotsc, h}$. By the assumption of our recursion, we have two cases: 

3.1. $ \card{ \qr{J_g}{J_g \cap J_{\ti{g}} } } = 1 $: This means 
\bgl
\label{3.1}
  J_g \subseteq J_{\ti{g}}.
\egl
If $J'_g = J_g \cap J'_{h+1}$, we must have $J'_g = J_g \cap J'_{h+1} \subseteq J'_{\ti{g}} $ so that (\ref{cond}) still holds.

If $J'_g = J_g \neq J_g \cap J'_{h+1}$, we can deduce $ \card{ \qr{J_g}{J_g \cap J'_{h+1} } } = \infty $. By (\ref{3.1}), we can consider the inclusion $\qr{J_g}{J_g \cap J'_{h+1}} \into \qr{J_{\ti{g}}}{J_{\ti{g}} \cap J'_{h+1}}$ which implies $\card{ \qr{J_{\ti{g}}}{J_{\ti{g}} \cap J'_{h+1} } } = \infty$. Thus, $J'_{\ti{g}} = J_{\ti{g}}$ and condition (\ref{cond}) holds by our assumption on $\gekl{J_1, \dotsc, J_h}$.

3.2. $ \card{ \qr{J_g}{J_g \cap J_{\ti{g}} } } = \infty $: 

If $J'_g=J_g$, consider the projection $ \qr{J_g}{J_g \cap J'_{\ti{g}}} \onto \qr{J_g}{J_g \cap J_{\ti{g}}} $ ($J'_{\ti{g}} \subseteq J_{\ti{g}}$ by construction) and deduce $ \card{ \qr{J'_g}{J'_g \cap J'_{\ti{g}} } } = \card{ \qr{J_g}{J_g \cap J_{\ti{g}} } } = \infty $. Thus, condition (\ref{cond}) is valid.

Otherwise, $J'_g = J_g \cap J'_{h+1} \neq J_g $ and we have the exact sequence 
\bgloz
  \qr{J'_g}{J'_g \cap J_{\ti{g}}} \into \qr{J_g}{J'_g \cap J_{\ti{g}}} \onto J_g / J'_g, 
\egloz
where $\card{J_g / J'_g} = \card{ \qr{J_g}{J_g \cap J'_{h+1}} } < \infty$ since $J'_g = J_g \cap J'_{h+1} \neq J_g $, and $\card{\qr{J_g}{J'_g \cap J_{\ti{g}}}} = \infty$ since $ \card{ \qr{J_g}{J_g \cap J_{\ti{g}}} } = \infty $ and because there is the canonical projection 
$ \qr{J_g}{J'_g \cap J_{\ti{g}} } = \qr{J_g}{J_g \cap J'_{h+1} \cap J_{\ti{g}}} \onto \qr{J_g}{J_g \cap J_{\ti{g}}} $.

Thus, $\card{ \qr{J'_g}{J'_g \cap J'_{\ti{g}}} } = \card{ \qr{J'_g}{J'_g \cap J_{\ti{g}}} } = \infty$ proving our second claim.
\eproof

Therefore, the projection $p \in \gekl{p_i}$ can be written as
\bgl
\label{p2}
  \sum_j n_j \e{a_j+I_j}- \sum_{j'} \ti{n}_{j'} \e{\ti{a}_{j'} + \ti{I}_{j'}}
\egl
with finitely many $n_j, \ti{n}_{j'} \in \Zz \pos$ so that 
\bgl
\label{a0}
  \card{ \qr{J}{J \cap \ti{J}} } \in \gekl{1, \infty} \text{ for any } J, \ti{J} \in \gekl{I_j, \ti{I}_{j'}}.
\egl
Moreover, we can certainly arrange 
\bgloz
  \e{a_j+I_j} \neq \e{\ti{a}_{j'}+\ti{I}_{j'}} \text{ for any } j, j' \text{ with } n_j \neq 0, \ti{n}_{j'} \neq 0.
\egloz
Since the sum (\ref{p2}) is finite, there exists an ideal $I \in \gekl{I_j}$ which is maximal among the $\gekl{I_j}$ with respect to inclusion. This means that for all $j$, $I \subseteq I_j$ implies $I=I_j$. 

\blemma
\label{l2}
If $I$ is maximal among $\gekl{I_j}$, it is already maximal among the larger set of ideals $\gekl{I_j, \ti{I}_{j'}}$.
\elemma

\bproof
Let us assume that there exists $\ti{I} \in \ti{I}_{j'}$ with $I \subseteq \ti{I}$ and $I \neq \ti{I}$. Because of our assumption (\ref{a0}) (see Lemma \ref{l1}), we can conclude that $ \card{ \ti{I} / I } = \infty $. $p$ is a projection, hence positive, and $\fD_{\cF} [R]$ is commutative, so we get
\bgloz
  \sum_j n_j \e{a_j+I_j} \geq \e{\ti{a}+\ti{I}}. 
\egloz
Multiplying this inequality by $\e{\ti{a}+\ti{I}}$ yields
\bgl
\label{i1}
  \sum_j n_j \e{\ti{a}+\ti{I}} \cdot \e{a_j+I_j} \geq \e{\ti{a}+\ti{I}}
\egl
where we choose $\ti{a} \in R$ corresponding to $\ti{I}$ as in the sum (\ref{p2}).

In the sum on the left-hand side, we only have contributions of those indices $j$ with $ (\ti{a}+\ti{I}) \cap (a_j+I_j) \neq \emptyset $. For such $j$, the intersection $ (\ti{a}+\ti{I}) \cap (a_j+I_j) $ is of the form $c_j + (\ti{I} \cap I_j)$ by Lemma \ref{int}. 

Thus, Relation II.(ii) allows us to rewrite (\ref{i1}) as 
\bgl
\label{i2}
  \sum_j n_j \e{c_j+(\ti{I} \cap I_j)} \geq \e{\ti{a}+\ti{I}}.
\egl
But now, we know that for all $j$, $\ti{I} \nsubseteq I_j$ because otherwise, $I \subseteq \ti{I} \subseteq I_j$ for some $j$ would imply $I=\ti{I}$ since $I$ is maximal among the $\gekl{I_j}$. This contradicts our assumption $I \neq \ti{I}$. Therefore, by (\ref{a0}), we must have $ \card{ \qr{\ti{I}}{\ti{I} \cap I_j} } = \infty $ for all $j$. 

In this situation, Corollary (\ref{finun}) implies $\bigcup_j (c_j + (\ti{I} \cap I_j)) \subsetneq \ti{a}+\ti{I}$ where the union is taken over the indices $j$ which contribute to the left-hand side of (\ref{i2}), there are only finitely many of these.
  
Thus, we can choose $r \in (\ti{a}+\ti{I}) \setminus \bigcup_j (c_j + (\ti{I} \cap I_j))$.

The corresponding projection $\e{r+(\bigcap_j (\ti{I} \cap I_j))}$ satisfies $\e{r+(\bigcap_j (\ti{I} \cap I_j))} \leq \e{\ti{a}+\ti{I}}$ since $r+(\bigcap_j (\ti{I} \cap I_j)) \subseteq \ti{a}+\ti{I}$. But we also have $\e{r+(\bigcap_j (\ti{I} \cap I_j))} \perp \e{c_j+(\ti{I} \cap I_j)}$ for all $j$ because $r$ does not lie in $\bigcup_j (c_j + (\ti{I} \cap I_j))$ by our choice. 

Therefore, multiplying (\ref{i2}) by $\e{r+(\bigcap_j (\ti{I} \cap I_j))}$ yields $0 \geq \e{r+(\bigcap_j (\ti{I} \cap I_j))}$ which is a contradiction.
\eproof

Finally, we are able to prove the following

\blemma[Technical Lemma]
Assume that all $I \in \cF$ satisfy $I \nsubseteq Z$. Then, there exist finitely many pairwise orthogonal projections $p_i$ in $\fA_{\cF} [R]$ with $C^* \rukl{ \gekl{p_i} } = C^* \rukl{ \menge{\e{m+I}}{\beta_{(m,I)} \neq 0} }$ (see (\ref{sum})).

Moreover, we can find $b_i' \in R \reg$, $a_i' \in R$ with $ \e{a_i'+(b_i')} \leq p_i \fa i $.
\elemma

\bproof
Recall that we had 
\bgloz
  \Theta(x) = s_d^* (\sum_{(m,I)} \beta_{(m,I)} \e{m+I}) s_d \text{ (see (\ref{sum}))},
\egloz
where the sum is taken over finitely many pairs $ (m,I) \in R \times \cF' $ with 
\bgloz
  \cF' = \menge{ \bigcap_{i=1}^n (b_i \cdot I_i)}{b_i \in R \reg, I_i \in \cF \cup \gekl{R}}.
\egloz 

Later on, it will be useful that
\bgl
\label{d}
  m+I \subseteq (d) \fa (m,I) \text{ with } \beta_{(m,I)} \neq 0.
\egl

First of all, Lemma \ref{l0} gives us suitable projections $p_i$. Choose $p \in \gekl{p_i}$ and write $p$ as 
\bgl
\label{p3}
  \sum_j n_j \e{a_j+I_j}- \sum_{j'} \ti{n}_{j'} \e{\ti{a}_{j'} + \ti{I}_{j'}}
\egl
with finitely many $n_j, \ti{n}_{j'} \in \Zz \pos$ and
\bgl
\label{a1}
  \e{a_j+I_j} \neq \e{\ti{a}_{j'}+\ti{I}_{j'}} \text{ for any } j, j' \text{ with } n_j \neq 0, \ti{n}_{j'} \neq 0
\egl
as well as
\bgl
\label{a2}
  \card{ \qr{J}{J \cap \ti{J}} } \in \gekl{1, \infty} \text{ for any } J, \ti{J} \in \gekl{I_j, \ti{I}_{j'}}.
\egl
This is possible by Lemma \ref{l1}.

Moreover, by Lemma \ref{l2}, we can choose $I \in \gekl{I_j}$ maximal with respect to inclusion and then, $I$ will automatically be maximal among $\gekl{I_j,\ti{I}_{j'}}$. 

Now, we choose $a \in R$, $n \in \Zz \pos$ so that $ n \cdot \e{a+I} $ appears as a summand in (\ref{p3}). Multiplying $p$ with $\e{a+I}$ gives 
\bgl
\label{p4}
  \e{a+I} \cdot p = n \cdot \e{a+I} + \sum_k n_k \e{c_k+(I \cap I_k)} - \sum_l \ti{n}_l \e{\ti{c}_l + (I \cap \ti{I}_l)} 
\egl
for some $c_k, \ti{c}_l \in R$ and $n_k, \ti{n}_l \in \Zz \pos$. Here, we used Relation II.(ii) as well as Lemma \ref{int}. 

We must have $ \card{\qr{I}{I \cap I_k}} = \infty $, $ \card{\qr{I}{I \cap \ti{I}_l}} = \infty $ for all $k, l$ since for $ I' \in \gekl{I_j, \ti{I}_{j'}} $, either $I'=I$ or $ \card{ \qr{I}{I \cap I'} } = \infty $ by (\ref{a2}), and in the first case, the only possible contribution after multiplication with $\e{a+I}$ is $ n \cdot \e{a+I} $ by (\ref{a1}). 

Thus, Corollary \ref{finun} implies 
\bgloz
  \bigcup_k (c_k+(I \cap I_k)) \cup \bigcup_l (\ti{c}_l+(I \cap \ti{I}_l)) \subsetneq a+I
\egloz
as there are only finitely many summands in (\ref{p4}). 

Again, take $ r \in (a+I) \setminus (\bigcup_k (c_k+(I \cap I_k)) \cup \bigcup_l (\ti{c}_l+(I \cap \ti{I}_l)))$. As in the proof of Lemma \ref{l2}, we conclude that 
\bglnoz
  && \e{r + (\bigcap_k (I \cap I_k) \cap \bigcap_l (I \cap \ti{I}_l))} \leq \e{a+I} \\
  && \e{r + (\bigcap_k (I \cap I_k) \cap \bigcap_l (I \cap \ti{I}_l))} \perp \e{c_k+(I \cap I_k)} \fa k \\
  && \e{r + (\bigcap_k (I \cap I_k) \cap \bigcap_l (I \cap \ti{I}_l))} \perp \e{\ti{c}_l+(I \cap \ti{I}_l)} \fa l.
\eglnoz
Therefore, multiplying (\ref{p4}) with $\e{r + (\bigcap_k (I \cap I_k) \cap \bigcap_l (I \cap \ti{I}_l))}$ yields
\bglnoz
  && \e{r + (\bigcap_k (I \cap I_k) \cap \bigcap_l (I \cap \ti{I}_l))} \cdot p \\
  &=& \e{r + (\bigcap_k (I \cap I_k) \cap \bigcap_l (I \cap \ti{I}_l))} \cdot \e{a+I} \cdot p \\
  &=& n \cdot \e{r + (\bigcap_k (I \cap I_k) \cap \bigcap_l (I \cap \ti{I}_l))} \\
  &+& \e{r + (\bigcap_k (I \cap I_k) \cap \bigcap_l (I \cap \ti{I}_l))} 
  \cdot (\sum_k n_k \e{c_k+(I \cap I_k)} - \sum_l \ti{n}_l \e{\ti{c}_l + (I \cap \ti{I}_l)}) \\
  &=& n \cdot \e{r + (\bigcap_k (I \cap I_k) \cap \bigcap_l (I \cap \ti{I}_l))}.
\eglnoz

As the product must be a projection again because $\fD_{\cF} [R]$ is commutative, $n$ must be equal to $1$ and we have proven $\e{r + (\bigcap_k (I \cap I_k) \cap \bigcap_l (I \cap \ti{I}_l))} \leq p$. 

Now, as $I \nsubseteq Z$ for all $I \in \cF$, we can find $b' \in R \reg \cap (\bigcap_k (I \cap I_k) \cap \bigcap_l (I \cap \ti{I}_l))$. Thus, $\e{r+(b')} \leq \e{r + (\bigcap_k (I \cap I_k) \cap \bigcap_l (I \cap \ti{I}_l))} \leq p$.

If we do this for all projections $p_i$, we will get, for all $i$, $a'_i \in R$ and $b'_i \in R \reg$ with $\e{a'_i+(b'_i)} \leq p_i$ as desired.

\eproof

\subsubsection{Proof of the main result}

Finally, we are ready to prove our main result on the inner structure of ring C*-algebras. 

\btheo
\label{main}
Let $R$ be a commutative ring with unit and let $Z$ be the set of zero-divisors in $R$. Assume that
\bgln
  \label{t1}
  && \bigcap_{b \in R \reg} (b) = (0), \\
  \label{t2}
  && \text{for all ideals } I \subseteq Z \text{ in } R \text{ and for all } b \in R \reg, \card{ \qr{(b)}{(b) \cap I} } = \infty.
\egln
Moreover, let $\cF$ be a family of ideals in $R$. 

Then, the ring C*-algebra $ \fA_{\cF} [R] $ is purely infinite and simple if and only if all $I \in \cF$ satisfy $I \nsubseteq Z$.
\etheo

We have two implications to prove: 

\bprop
\label{if}
Assume that in the situation of the theorem, all $I \in \cF$ satisfy $I \nsubseteq Z $. Then $\fA_{\cF} [R]$ is purely infinite and simple.
\eprop

\bproof
We will establish the criterion of Proposition \ref{cri} for $A = \fA_{\cF} [R]$, $\cA = \cA_{\cF} [R] $ (see Lemma \ref{dense}) and the faithful conditional expectation $\Theta$ of Lemma \ref{fce}. 

So, let $x$ be a nontrivial element in $\cA_{\cF} [R]_{sa}$. Write 
\bgl
\label{x}
  x = \sum_{k,k',l,l',J} \alpha_{(k,k',l,l',J)} s_l^* u^k \e{J} u^{-k'} s_{l'}.
\egl
Just as in our Technical Lemma, we can write 
\bgloz
  \Theta (x) = s_d^* (\sum_{(m,I)} \beta_{(m,I)} \e{m+I}) s_d
\egloz
with $(m,I) \in R \reg \times \cF'$ and $m+I \subseteq (d)$ for all $(m,I)$ with $\beta_{(m,I)} \neq 0$. 

By the Technical Lemma, we can find finitely many pairwise orthogonal (nontrivial) projections $\gekl{p_i}$ with  $C^* \rukl{ \menge{\e{m+I}}{\beta_{(m,I)} \neq 0} } = C^* \rukl{ \gekl{p_i} }$. Furthermore, there exist $a'_i \in R$, $b'_i \in R \reg$ such that $ \e{a'_i+(b'_i)} \leq p_i $ for all $i$. Now, since $ \e{a'_i+(b'_i)} \leq p_i \leq \e{(d)} $, we must have $ a'_i + (b'_i) \subseteq (d) $ as we can look at the corresponding projections on $\ell^2(R)$ using the canonical projection $\pi: \fA_{\cF} [R] \lori \fA_{r,\cF} [R] \subseteq \cL (\ell^2(R))$, and in $\fA_{r,\cF} [R]$, $\E{a'_i+(b'_i)} \leq \E{(d)}$ is equivalent to $a'_i+(b'_i) \subseteq (d)$. 

Thus, the projections $F_i \defeq s_d^* \e{a'_i+(b'_i)} s_d = \e{a'_i/d+(b'_i/d)}$ satisfy $F_i \leq s_d^* p_i s_d$ and  
\bgl
\label{iso}
  C^* \rukl{\menge{\e{m+I}}{\beta_{(m,I)} \neq 0}} = C^* \rukl{\gekl{p_i}} \lori C^* \rukl{\gekl{F_i}}; y \loma \sum_i F_i y F_i
\egl
is an isomorphism as $s_d^* p_i s_d$ is mapped to $F_i$. 

These projections $F_i$ have all the desired properties except that 
\bgloz
  \sum_i F_i x F_i = \sum_i F_i \Theta(x) F_i
\egloz
does not hold in general. 

So, we have to look at the representation (\ref{x}) of $x$ again:

Call a multiindex $(k,k',l,l',J)$ critical if $ \alpha_{(k,k',l,l',J)} s_l^* u^k \e{J} u^{-k'} s_{l'} \neq 0 $ and $ \delta_{k,k'} \delta_{l,l'} = 0 $. There are finitely many of them, but these are the summands of $x$ in (\ref{x}) which are of interest since 
\bgl
\label{diff}
  x-\Theta(x) = \sum_{(k,k',l,l',J) \text{ critical}} \alpha_{(k,k',l,l',J)} s_l^* u^k \e{J} u^{-k'} s_{l'}.
\egl

We claim: For each $i$ there exists $a_i \in a'_i/d + (b'_i/d)$ with $k-k'-(l-l')a_i \neq 0$ for all critical indices $(k,k',l,l',J)$. 

Proof of the claim:

For $l=l'$ we must have $k \neq k'$ as $ \delta_{k,k'} \delta_{l,l'} = 0 $. Thus, $k-k'-(l-l')a_i$ will be nontrivial for any $a_i$. So, we only have to consider critical indices with $l \neq l'$. 

Assume that our assertion is false, this amounts to saying: 

There exists $i$ with 
\bgl
\label{co1}
  a'_i/d+(b'_i/d) \subseteq \bigcup_{(k,k',l,l',J) \text{ critical with } l \neq l'} A(k,k',l,l')
\egl
where $ A(k,k',l,l') = \menge{r \in R}{k-k'-(l-l')r = 0} $.

Obviously, we only have to consider those critical multiindices $(k,k',l,l',J)$ with $l \neq l'$ and $A(k,k',l,l') \neq \emptyset$, let us call them hypercritical. Now, for any hypercritical multiindex $(k,k',l,l',J)$, take some $r(k,k',l,l') \in A(k,k',l,l')$. It is clear that 
\bgl
\label{A}
  A(k,k',l,l') = r(k,k',l,l') + (0:(l-l'))
\egl
where $(0:(l-l')) = \menge{r \in R}{(l-l')r=0}$. 

Substituting (\ref{A}) into (\ref{co1}), we get 
\bgl
\label{co2}
  a'_i/d+(b'_i/d) \subseteq \bigcup_{(k,k',l,l',J) \text{ hypercritical}} (r(k,k',l,l')+(0:(l-l'))).
\egl
Since $b'_i/d \in R \reg$ and $(0:(l-l')) \subseteq Z$, we have 
\bgloz
  \card{ \qr{(b'_i/d)}{(b'_i/d) \cap (0:(l-l'))} } = \infty
\egloz
by the second condition (\ref{t2}) on our ring $R$. 

But then, (\ref{co2}) contradicts Corollary \ref{finun}. This proves our claim. 

Now, we have seen that it is possible to choose for each $i$ some $a_i \in a'_i/d+(b'_i/d)$ with $k-k'-(l-l')a_i \neq 0$ for all critical multiindices $(k,k',l,l',J)$. Because of the condition (\ref{t1}) on $R$, we can find, for each critical multiindex and for every $i$, an element $b(k,k',l,l',J,i) \in R \reg$ with 
$b(k,k',l,l',J,i) \ntei (k-k'-(l-l')a_i)$. 

Then, define for all $i$
\bgl
\label{b_i}
  b_i \defeq b'_i/d \cdot \prod_{(k,k',l,l',J) \text{ critical}} b(k,k',l,l',J,i) \in R \reg.
\egl

Our final claim is that the projections $f_i \defeq \e{a_i+(b_i)} $ satisfy the criterion of Proposition \ref{cri}.

Proof of the final claim: Recall that we have to show
\bgln
\label{1}
  && f_i \perp f_j \fa i \neq j; f_i \overset{s_i}{\sim} 1 \text{ via isometries } s_i \in A \fa i, \\
\label{2}
  && \norm{\sum_i f_i \Theta(x) f_i } = \norm{\Theta(x)}; f_i x f_i = f_i \Theta(x) f_i \in \Cz f_i \fa i.  
\egln

To show (\ref{1}), just note that for all $i$, $f_i \leq F_i$ by construction. This implies $f_i \perp f_j \fa i \neq j$. 
Moreover, $f_i = (u^{a_i} s_{b_i}) (u^{a_i} s_{b_i})^* \sim 1$ for all $i$. 

To show the first part of (\ref{2}), we deduce from $f_i \leq F_i$ for all $i$ that 
\bgloz
  C^* \rukl{ \menge{\e{m+I}}{\beta_{(m,I)} \neq 0} } \lori C^* \rukl{ \gekl{f_i} }; y \loma \sum_i f_i y f_i 
\egloz
is an isomorphism because (\ref{iso}) is an isomorphism. Hence, this map is isometric and we get $\norm{\sum_i f_i \Theta(x) f_i } = \norm{ \Theta(x) }$ as well as $f_i \Theta(x) f_i \in \Cz f_i$ for all $i$. 

It remains to show $f_i x f_i = f_i \Theta(x) f_i$ for all $i$. Using (\ref{diff}), we get
\bglnoz
  && f_i(x-\Theta(x))f_i \\
  &=& f_i (\sum_{(k,k',l,l',J) \text{ critical}} \alpha_{(k,k',l,l',J)} s_l^* u^k \e{J} u^{-k'} s_{l'}) f_i \\
  &=& \sum \alpha_{(k,k',l,l',J)} s_l^* u^k (u^{-k} s_l f_i s_l^* u^k) \e{J} (u^{-k'} s_{l'} f_i s_{l'}^* u^{k'}) u^{-k'} s_{l'} \\
  &=& \sum \alpha_{(k,k',l,l',J)} s_l^* u^k \e{-k+la_i+(lb_i)} \cdot \e{-k'+l'a_i+(l'b_i)} \cdot \e{J} u^{-k'} s_{l'}.
\eglnoz
Now, if $ \eckl{-k+la_i+(lb_i)} \cap \eckl{-k'+l'a_i+(l'b_i)} \neq \emptyset $, we are able to deduce that $-k'+l'a_i-(-k+la_i) = k-k'-(l-l')a_i \text{ lies in } (b_i)$ contradicting our choice of $b_i$ (see (\ref{b_i})). Thus, $ f_i(x-\Theta(x))f_i = 0 $ for all $i$. This proves the final claim. 

Now, our Proposition follows from Proposition \ref{cri}.

\eproof

This has certainly been the main part of the proof of the Theorem \ref{main}. It remains to show the \an{only if}-part.

\bprop
Let $R$ be a commutative ring with unit satisfying the condition
\bgl
\label{r}
  \text{for all ideals } I \subseteq Z \text{ in } R \text{ and for all } b \in R \reg, \card{ \qr{(b)}{(b) \cap I} } = \infty
\egl
where $Z$ is the set of zero-divisors in $R$. Let $\cF$ be a family of ideals in $R$ and assume that there exists $J \in \cF$ with $J \subseteq Z$. 

Then, the ideal $\fI$ generated by $\e{J}$ in $\fA_{\cF} [R]$ is a proper ideal.
\eprop

\bproof
We certainly have $(0) \neq \fI $, so it remains to show $ \fA_{\cF} [R] \neq \fI $. 

First, take $\cI \defeq \lspan \rukl{ \menge{s_l^* u^{-k} \e{\ti{J}} u^{k'} s_{l'}}{k, k' \in R; l, l' \in R \reg; \ti{J} \in \cF', \ti{J} \subseteq J} }$, 
where $\cF' = \menge{\bigcap (b_i J_i)}{b_i \in R \reg, J_i \in \cF \cup \gekl{R}}$. 

We claim $ \overline{\cI} = \fI $. Proof of the claim:

As $\e{\ti{J}} = \e{J} \cdot \e{\ti{J}} \in \fI$, we have $\cI \subseteq \fI$ and thus $\overline{\cI} \subseteq \fI$. For the other inclusion, we show that $\cI$ is a two-sided ideal of $\cA_{\cF} [R]$ in the algebraic sence. This then implies $\overline{\cI} \supseteq \fI$ since $ \overline{\cA_{\cF} [R]} = \fA_{\cF} [R] $ (by Lemma \ref{dense}), and thus, $\overline{\cI}$ is a two-sided ideal in $\fA_{\cF} [R]$ with $\overline{\cI} \supseteq \cI \ni \e{J}$. 

To see that $\cI \cdot \cA_{\cF} [R] \subseteq \cI$, take $s_b^* u^a \e{\ti{J}} u^{-a'} s_{b'} \in \cI, s_d^* u^c \e{I} u^{-c'} s_{d'} \in \cA_{\cF} [R]$, and by the computation of Lemma \ref{dense}, their product is of the form
\bgloz
  s_{bd}^* u^{(a-a')d+\ti{a}} \e{ (d \cdot \ti{J}) \cap (db') \cap (b' \cdot I) } u^{-\ti{a}+(c-c')b'} s_{b'd'}
\egloz
which lies in $\cI$ since $ (d \cdot \ti{J}) \cap (db') \cap (b' \cdot I) \subseteq \ti{J} \subseteq J $. Thus, $\overline{\cI}=\fI$. 

Now, we pass over to $\fA_{r,\cF} [R]$ and look at $\pi(\fI) \subseteq \cL(\ell^2(R))$. As $\overline{\cI} = \fI$, we must have $\overline{ \pi(\cI) }= \pi(\fI)$. Now, an arbitrary element of $\pi(\cI)$ is a finite sum of the form $\sum \alpha_{(a,a',b,b',\ti{J})} S_b^* U^{-a} \E{\ti{J}} U^{a'} S_{b'} \text{ with } \ti{J} \subseteq J \subseteq Z$. Its support projection is dominated by $\E{\bigcup ((a'+\ti{J}):b')}$ where $((a'+\ti{J}):b') = \menge{r \in R}{b'r \in a'+\ti{J}}$. Note that the unions we are taking are always finite. The sets $ ((a'+\ti{J}):b') $ are either empty or of the form $r+(\ti{J}:b')$ for some $r \in ((a'+\ti{J}):b')$. Thus, the support projection of any element in $\pi(\cI)$ is dominated by a projection of the form $\E{\bigcup_{i=1}^n (r_i + \ti{J}_i)}$ with $\ti{J}_i \subseteq J \subseteq Z$. 

But this projection cannot be $1$ as $\bigcup_{i=1}^n (r_i + \ti{J}_i) \neq R$ by (\ref{r}) and Corollary \ref{finun}. Therefore, $\pi(\cI)$ does not contain any invertible element, and thus, $1 \notin \pi(\fI) = \overline{\pi(\cI)}$. Hence it follows that $1 \notin \fI$ in $\fA_{\cF} [R]$ and thus $\fI \neq \fA_{\cF} [R]$.
\eproof

Let us add some remarks. 

\bcor
If $R$ is an integral domain, then $\fA_{\cF} [R]$ is purely infinite and simple if and only if $R$ is not a field and $(0) \notin \cF$.
\ecor
\bproof
We just have to translate conditions (\ref{t1}) and (\ref{t2}) of the Theorem to this special case. 

First of all, note that if $R$ is a field, $\cF$ will be $\emptyset$ or $\gekl{R}$ and $\fA_{\cF} [R]$ will then coincide with the full group C*-algebra $C^*(P_R)$ which cannot be simple. 

If $R$ is not a field, (\ref{t1}) is automatically satisfied as any nonzero element in $R$ is not divisible by itself times a noninvertible element. Moreover, $R$ must be infinite and thus, we always have $\# (b) = \infty$ for all $b \in R \reg$. But as $Z=\gekl{0}$, this shows that $R$ satisfies (\ref{t2}) as well. Furthermore, as $Z=\gekl{0}$, the statements \an{$(0) \notin \cF$} and \an{for all $I \in \cF$: $I \nsubseteq Z$} are equivalent.
\eproof

Moreover, note that condition (\ref{t2}) is always satisfied if $R$ contains an infinite field $K$ as then, any quotient carries the structure of a $K$-vectorspace.

We can also deduce

\bcor
\label{A-Ar}
In the situation of Theorem \ref{main} and under the assumption that all $I \in \cF$ satisfy $I \nsubseteq Z$, the canonical epimorphism $\pi: \fA_{\cF} [R] \lori \fA_{r, \cF} [R]$ is an isomorphism. 

In particular ($\cF = \emptyset$), if $R$ satisfies the conditions (\ref{t1}) and (\ref{t2}) as in Theorem \ref{main}, $\pi: \fA [R] \lori \fA_r [R]$ is always an isomorphism. 
\ecor

\section{Examples}
\label{E}

As we have explained, the first examples we looked at were the rings of integers in number fields. This type of examples has already been discussed in Section \ref{Re}. Now, let us discuss three examples closely related to rings of integers, namely their localizations, group rings and rings of matrices. These are rings with finite quotients. Furthermore, we will also look at rings with infinite quotients, namely $\Qz [T]$ and $\Zz[i \sqrt{5}] [T]$. 

The first three examples have in common that once we go over to the minimal automorphic dilation, our construction behaves very naturally: Inverting a prime corresponds to leaving out the corresponding place in the finite adele ring, and taking group rings or matrices corresponds to taking group rings or matrices with coefficients in the finite adele ring. Moreover, for $\Qz [T]$, our construction yields a natural generalization of the finite adele ring sharing many structural properties with the classical one for global fields but lacking a reasonable ring structure. And our investigation of $\Zz [i \sqrt{5}] [T]$ reveals an interesting phenomenon involving the conditional expectation of Lemma \ref{fce}. This is due to the fact that $\Zz [i \sqrt{5}]$ has class number greater than $1$. 

\subsection{$\Zz [\tfrac{1}{p}]$}

Let us take for instance $\Zz$ and formally invert the prime $p$. In the corresponding ring C*-algebra of $R = \Zz [\tfrac{1}{p}]$, we have made $s_p$ unitary, which means $e_p = 1$. Thus, $\Spec \fD [R] \cong \prod_{q \neq p} \Zz_q$ and by Lemma \ref{cor}, 
\bgloz
  \fA [\Zz [\tfrac{1}{p}]] \cong C(\prod_{q \neq p} \Zz_q) \rte \Zz [\tfrac{1}{p}] \rtimes (\Zz [\tfrac{1}{p}]) \reg 
  \sim_M C_0({\prod_{q \neq p}}' \Qz_q) \rtimes \Qz \rtimes \Qz \reg.
\egloz

\subsection{$\Zz [\Zz / p \Zz]$}

Let $p$ be a prime and $t$ the generator of $\Zz / p \Zz$. The group ring $R = \Zz [\Zz / p \Zz]$ is commutative with zero-divisors $Z = (1+\dotsb+t^{p-1}) \cup (1-t)$. $R$ has the property that for each $b \in R \reg$, there exists $d \in R \reg$ with $d \cdot b \in \Zz \pos \cdot 1$. Thus, 
\bgloz
  \Spec \fD [R] \cong \plim_{m \in \Zz \pos} R / (m) \cong \plim_{m \in \Zz \pos} (\Zz / m \Zz) [\Zz / p \Zz] \cong \widehat{\Zz} [\Zz / p \Zz] 
\egloz
and $\fA [ \Zz [\Zz / p \Zz] ] \cong C(\widehat{\Zz} [\Zz / p \Zz]) \rte (\Zz [\Zz / p \Zz]) \rtimes (\Zz [\Zz / p \Zz]) \reg$. The minimal automorphis dilation is given by the action of $(\Qz [\Zz / p \Zz]) \rtimes (\Qz [\Zz / p \Zz]) \reg$ on $\Az_f [\Zz / p \Zz]$ via affine transformations. As $\Qz [\Zz / p \Zz] \cong \Qz \times \Qz [\zeta]$ with a primitive $p$-th root of unity $\zeta$, we get by Lemma \ref{cor} 
\bglnoz
  && \fA[M_l (\so)] \sim_M C_0(\Az_f[\Zz / p \Zz]) \rtimes (\Qz [\Zz / p \Zz]) \rtimes (\Qz [\Zz / p \Zz]) \reg \\
  &\cong& (C_0(\Az_{f, \Qz}) \rtimes \Qz \rtimes \Qz \reg) \otimes (C_0(\Az_{f, \Qz [\zeta]}) \rtimes \Qz [\zeta] \rtimes \Qz [\zeta] \reg).
\eglnoz

\subsection{$M_k (\so)$}

Let $\so$ be the ring of integers in a number field $K$. We would like to look at the ring of $l \times l$-matrices over $\so$. In this  mildly noncommutative situation, our construction still produces a natural C*-algebra. 

Let $R = M_l (\so)$. We have $R \reg = \menge{b \in M_l (\so)}{\det (b) \neq 0}$. Again, for each $b \in R \reg$, there exists $d \in R \reg$ with $bd = db \in \Zz \pos \cdot 1_l \subseteq M_l (\so)$. Thus, 
\bgloz
  \Spec \fD[R] \cong \plim_{b \in R \reg} R/(b) \cong \plim_{m \in \Zz \pos} R/(m) \cong \plim_{m \in \Zz \pos} M_l(\so/(m)) \cong M_l(\widehat{\so}).
\egloz

So, all in all, we get $\fA [ M_l (\so) ] \cong C(M_l (\widehat{\so})) \rte R \rtimes R \reg$. The associated minimal automorphic dilation is given by $C_0(M_l(\Az_f))$ together with the action of $M_l(K) \rtimes GL_l(K)$ where $M_l(K)$ acts additively and $GL_l(K)$ acts via matrix multiplication (from the left). Thus, by Lemma \ref{cor},
\bgloz
  \fA[M_l(\so)] \sim_M C_0(M_l(\Az_f)) \rtimes M_l(K) \rtimes GL_l(K).
\egloz

\subsection{The polynomial ring with rational coefficients}

Consider the ring $ R = \Qz [T] $. This example has the following nice properties:

1. $R$ contains an infinite field, namely $ \Qz $. This implies that for any ideals $ I \subseteq J $ in $R$, we must have $ \card{ J/I } \in \lge 1, \infty \rge $ because we are considering vector spaces over $ \Qz $. Thus, most of the manipulations of Section \ref{mani} can be simplified.

2. $R$ is a principal ideal domain. As a consequence, any choice of the family $ \cF $ must lead to the C*-algebra $\fA [R]$ corresponding to the choice $\cF = \emptyset$. In addition, this property of $R$ ensures that the C*-subalgebras $ \fD [R] $ and $ \Theta \lru \fA [R] \rru $ must coincide because 
\bgloz
  s_d^* \e{ a+(b) } s_d = s_d^* \e{ (d) \cap (a+(b)) } s_d = s_d^* \e{ a'+(b') } s_d = \e{ a'/d + (b'/d) } \in \fD [R] 
\egloz
for some $a', b' \in (d)$ because of Lemma \ref{int} and the fact that $R$ is a principal ideal domain. At this point, we should probably point out that if one thinks of $\Theta(\fA[R])$ as the fixed point algebra with respect to the action of $R \rtimes R \reg$, this statement is not surprising as we have $\fA[R] \cong \fD[R] \rte (R \rtimes R \reg)$. But we are in the case of semigroups now, and an equality like $ \fD [R] = \Theta \lru \fA [R] \rru $ (which is well-known for group actions) will only be valid in nice cases, as we will see.

Now, let us describe the spectrum of the commutative C*-subalgebra $ \fD [R] $. First, choose a set of (pairwise nonassociated) representatives $ \lge p_i \rge $ of irreducible polynomials in $R$.

\bprop
\label{SpecQT}
$ \Spec \fD [R] $ can be identified with $\plim_n \gekl{Z_n ; \pi_{n+1,n}}$ where we have 
\bgloz
  Z_n = \underset{0 \leq i_1, \dotsc,i_n \leq n}{\dotcup} R/(p_1^{i_1} \dotsm p_n^{i_n})
\egloz 
as sets. Moreover, the structure maps $\pi_{n+1,n}$ are given by the projections
\bgloz
  R / \lru p_1^{i_1} \dotsm p_n^{i_n} p_{n+1}^{i_{n+1}} \rru \lori R / \lru p_1^{\min(i_1,n)} \dotsm p_n^{\min(i_n,n)} \rru.
\egloz
The topology of $ Z_n $ can be described as follows: 

For any $ z \in Z_n $ define $ I_z $ to be the ideal of $R$ such that $ z \in R / I_z \subseteq Z_n $. We say that $ z \leq z' $ if $ I_z \subseteq I_{z'} $ and if 
$ z + I_{z'} = z' + I_{z'} $. This defines an order relation on $ Z_n $. Now, the topology on $ Z_n $ is given by the property that a sequence $ (x_m) $ converges to 
$ x \in R/I $ in $Z_n$ if and only if $x$ is the only minimal element of $\menge{x' \in Z_n}{x_m \leq x' \text{ for almost all } m}$ with respect to $\leq$.
\eprop

\bproof
We can write $ \fD [R] \cong \ilim_n \fD_n $ with
\bgloz
  \fD_n \defeq C^* \lru \lge \e{ a + \lru p_1^{i_1} \dotsm p_n^{i_n} \rru } \dop 0 \leq i_1, \dotsc, i_n \leq n \rge \rru.
\egloz
Thus, $\Spec \fD [R] \cong \plim_n \Spec \fD_n$. 

Now, associate to $z \in Z_n$ the projection $e_z \defeq \e{z+I_z} \in \fD_n$. Every $\chi \in \Spec \fD_n$ corresponds bijectively to an element $z(\chi) \in Z_n$ with 
\bgloz
  \chi ( e_{z'} ) = 
  \bfa 
    1 \falls z' \geq z( \chi ) \\
    0 \sonst.
  \efa
\egloz
This gives the bijection between $ \Spec \fD_n $ and $ Z_n $, denoted by $ \chi \mapsto z( \chi ) $. Let us denote the inverse by $ \chi_z \mapsfrom z $.

Furthermore, a sequence $ (\chi_{x_m}) $ converges to $ \chi_x $ if and only if 
\bgloz
  \lim_{m \lori \infty} \chi_{x_m} (e_{x'}) = \chi_x ( e_{x'} ) \fa x' \in Z_n. 
\egloz
But this is equivalent to the statement \an{$x' \geq x \LRarr x' \geq x_m$ for almost all $m$} which is the same as saying that $x$ is the only minimal element in $\menge{x' \in Z_n}{x_m \leq x' \text{ for almost all } m}$. Thus, $ \Spec \fD_n $ and $ Z_n $ are homeomorphic. 

The description for the structure maps can be instantly deduced from these identifications.
\eproof

Moreover, we even have the following
\bprop
There is an embedding $R \into \Spec \fD [R]$ with dense image.
\eprop

\bproof
Let $ I_n $ be the ideal $(p_1^n \dotsm p_n^n)$ of $R$. 
The embedding is given as the composition 
\bgloz
  \iota : R \into \plim_n R / I_n \into \Spec \fD [R]. 
\egloz
To see that $ \iota $ has dense image, take any $z$ in $\Spec \fD [R] \cong \plim_n \gekl{Z_n;\pi_{n+1,n}}$, which we view as an element $(z_n) \in \prod_n Z_n$ with $\pi_{n+1,n} (z_{n+1}) = z_n$ for all $n \in \Zz \pos$. Now, for all $ n \in \Zz_{>0} $ we can find a sequence $ z_n^{(m)} \in R $ such that 
\bgloz
  \lim_{m \rarr \infty} p_n \lru z_n^{(m)} \rru = z_n \in Z_n 
\egloz
where $p_n$ is the canonical projection $ R \ri R / (p_1^n \dotsm p_n^n)$. Thus, we can choose a suitable diagonal subsequence $(z_n^{(m_n)})_n$ in $R$ with $\lim_{n \ri \infty} z_n^{(m_n)} = z$.
\eproof

Moreover, we can define $\fD_k = C^*(\menge{\e{a+(p_k^i)}}{a \in R, i \in \Zz_{>0}})$, so that 
\bgloz
  \fD [R] \cong \ilim_N \otimes_{k=1}^N \fD_k \text{ and }
  \Spec \fD [R] \cong \plim_N \prod_{k=1}^N \Spec \fD_k \cong \prod_{k=1}^{\infty} \Spec \fD_k.
\egloz
Essentially the same argument used in the proof of Proposition \ref{SpecQT} shows that $\Spec \fD_k \cong \plim_M \dotcup_{i=0}^M R / (p_k^i)$ with a similar description of the topology as in Proposition \ref{SpecQT}.

But in contrast to all these similarities, there is also one striking feature of the infinite case which did not occur before: There is no ring structure on $ \Spec \fD [R] $ extending the one of $R$ and compatible with the topology on the spectrum. The reason is that we have added several \an{points at infinity} constructing $ \Spec \fD [R] $ out of the quotients $ R / I_n $. We were forced to do so in order to compactify $ \plim_n R / I_n $ (which is not even locally compact). In the finite case, this problem does not occur.

\subsection{The polynomial ring with coefficients in $ \Zz [i\sqrt{5}] $}

Our last example is the ring $R = \so [T]$ where $\so = \Zz [i \sqrt{5}]$ is the ring of integers in $\Qz [i \sqrt{5}]$. This ring is of interest as $\so$ is not a unique factorization domain. Equivalently, we can also say that the class number of $\Qz [i \sqrt{5}]$ is strictly greater than $1$. This leads to a new phenomenon which we did not encounter in our first example.

\bprop
The commutative C*-subalgebras $ \fD [R] $ and $ \Theta ( \fA [R] ) $ do not coincide.
\eprop

\bproof
The idea is to use the two factorizations $2 \cdot 3 = (1+i\sqrt{5}) \cdot (1-i\sqrt{5})$ to construct a projection in $\Theta(\fA [R])$ which does not lie in $\fD [R]$. 

Take $p \defeq s_2^* \e{(1+i\sqrt{5})} s_2$. By property (\ref{condexp}) of $\Theta$, $p$ lies in $\in \Theta(\fA [R])$. 

It remains to show that $p$ does not lie in 
\bgloz
  \fD [R] = C^*(\menge{\e{a+(b)}}{a \in R, b \in R \reg}) = \clspan(\menge{\e{a+\bigcap_{i=1}^n (b_i)}}{a \in R, b_i \in R \reg}).
\egloz
Assume the contrary. Using functional calculus, we obtain a sequence of projections in $\lspan (\menge{\e{a + \bigcap_{i=1}^n (b_i)}}{a \in R, b_i \in R \reg})$ converging to $p$. Hence, $p$ itself must already lie in $\lspan (\menge{\e{a + \bigcap_{i=1}^n (b_i)}}{a \in R, b_i \in R \reg})$.

Thus, by similar arguments as in Section \ref{mani}, we know that $p$ can be written as a (finite) linear combination of the projections $ \e{ a + \bigcap_j (b_j) } $ with integer coeffients. Moreover, if the projection $1$ appears with nontrivial coeffient, it must have coefficient $1$ (see the proof of the Technical Lemma). 
Therefore, 
\bgloz
  p \leq N \sum_{j=1}^n \e{ a_j + (b_j) } \text{ or } 1-p \leq N \sum_{j=1}^n \e{a_j + (b_j)}
\egloz
for some $N \in \Zz \pos$, $a_j \in R$, $b_j \in R \reg \setminus R^*$. Now, we deduce from this with the help of the regular representation of $\fA [R]$ on $\ell^2 (R)$ that 
\bgl
\label{inc}
  I \subseteq \bigcup_{j=1}^n (a_j + (b_j)) \text{ or } I^{\complement} \subseteq \bigcup_{j=1}^n (a_j + (b_j))
\egl
in $ R = \so [T] $, where $I \defeq [(1+i\sqrt{5}):2] = \menge{r \in R}{2r \in (1+i\sqrt{5})}$ is the ideal of $R$ which corresponds to $p$ in $ \fA_r [R] $. But from the second inclusion of (\ref{inc}), we can deduce 
\bgloz
  1 + I \subseteq I^{\complement} \subseteq \bigcup_{j=1}^n (a_j + (b_j))
\egloz 
which implies $I \subseteq \bigcup_{j=1}^n ((a_j-1) + (b_j))$. Therefore, we can assume without loss of generality that $I \subseteq \bigcup_{j=1}^n (a_j + (b_j))$ holds for some $ a_j \in R $, $ b_j \in R \reg \setminus R^* $.

Now, if we can show that $\card{I / (I \cap (b_j))} = \infty \fa j$, we will get the desired contradiction by Lemma \ref{finun}.

To this end, define $ \pi : R \lori \so $ as the canonical projection ($ \pi = \ev_0 $) and $ \iota $ as the inclusion $ \so \into R $. 
We have the following cases to consider:

1. $ b_j \in \so $. We know that $ \pi ( I ) \nsubseteq (b_j) $ in $ \so $ because otherwise, $ 1+i\sqrt{5} \in \pi (I) $ would imply $ (b_j) = (1+i\sqrt{5}) $ 
(because $ 1+i\sqrt{5} $ is irreducible) which contradicts $ 3 \in \pi (I) $, $ 3 \notin (1+i\sqrt{5}) $. This shows that we can find $r \in \pi(I)$ which does not lie in $(b_j)$. Now, consider the infinite subset $\menge{r \cdot T^i}{i \in \Zz \pos}$ of $I$, we have for all $k \neq l$ in $\Zz \pos$: $r \cdot T^k - r \cdot T^l \notin (b_j)$ because a polynomial lies in $(b_j) \subseteq R$ if and only if its coefficients lie in $(b_j) \subseteq \so$. 

Thus, we can conclude that $\card{ I / (I \cap (b_j))} = \infty$.

2. $ b_j \notin \so $. Then, $ \de (b_j) \geq 1 $. Hence we have $ \iota ( \pi (I) ) \cap (b_j) = (0) $, which means that any nontrivial element of $ \iota ( \pi (I) ) $ lies in $I$, but not in $ (b_j) $. Therefore, we can conclude that $ \card{ \lru I + (b_j) \rru / (b_j) } = \infty $. But this means $\card{ I / (I \cap (b_j))} = \infty$ as we always have $ I / \lru I \cap (b_j) \rru \cong \lru I + (b_j) \rru / (b_j) $.

So, all in all, we have shown that $p = s_2^* \e{(1+i\sqrt{5})} s_2 \in \Theta(\fA[R]) \setminus \fD [R]$.
\eproof

As a last comment, we remark that this argument can be refined to show that in any number field $K$, $\fD[\so_K [T]] = \Theta(\fA[\so_K [T]])$ holds if and only if 
$h_K = 1$. 

\appendix
\section{Semigroup C*-algebras and semigroup crossed products}

We would like to recall some basic definitions concerning (discrete) semigroup crossed products by endomorphisms. Moreover, we explain how one can take up the ideas of considering certain subsets as underlying data (as in Section \ref{R}) to generalize the notion of semigroup C*-algebras. We can even go further and consider semigroup crossed products by automorphisms. 

\subsection{On crossed products by discrete semigroups}
\label{semicropro}

Here are just some basic definitions on crossed products by discrete semigroups.

First of all, a semigroup is a set $P$ together with an associative binary operation 
\bgloz
  P \times P \lori P; (p,q) \loma p q.
\egloz

A unit is an element $e \in P$ with $e p = p e = p$ for all $p \in P$. If such a unit exists, it is unique. $P$ is called left-cancellative if for any $p, q, q' \in P$, $pq = pq'$ implies $q=q'$. A semigrouphomomorphism is a map between two semigroups which respects the binary operations.

From now on, we will only consider semigroups with unit. All semigrouphomomorphisms shall respect the unit elements. Moreover, we will not be concerned with topologies on $P$, which means that we are talking about the discrete case. 

\bdefin
A C*-dynamical semisystem is a triple $(D,P,\alpha)$ consisting of a unital C*-algebra $D$, a semigroup $P$ and a semigrouphomomorphism $\alpha: P \lori \End(D)$.
\edefin

\bdefin
A covariant representation of a C*-dynamical semisystem $(D,P,\alpha)$ is a triple $(A,\pi,\rho)$ consisting of a unital C*-algebra $A$, a unital homomorphism $ \pi: D \lori A$ and a homomorphism of semigroups $\rho: P \lori \Iso(A)$ such that the following covariance relation is fulfilled:

For all $d \in D$, $p \in P$, the equation $\rho(p) \pi(d) \rho(p)^* = \pi(\alpha(p)(d))$ holds.
\edefin

\bdefin
\label{mor}
A morphism $\Phi \dop (A_1, \pi_1, \rho_1) \lori (A_2, \pi_2, \rho_2)$ of two covariant representations $(A_i, \pi_i, \rho_i)$ ($i = 1, 2$) of $(D,P,\alpha)$ is a (necessarily unital) homomorphism $A_1 \overset{\Phi}{\lori} A_2 \text{ such that } \Phi \circ \pi_1 = \pi_2 \text{ and } \Phi \circ \rho_1 = \rho_2$.
\edefin

\bdefin
\label{defcropro}
Let $(D,P,\alpha)$ be a C*-dynamical semisystem. The crossed product associated to this C*-dynamical semisystem is the covariant representation 
$(D \rte_\alpha P,j,w)$ of $(D,P,\alpha)$ satisfying the following universal property:

For any covariant representation $(A,\pi,\rho)$ of $(D,P,\alpha)$, there exists a uniquely determined morphism of covariant representations 
\bgloz
  \Phi_{(A,\pi,\rho)}: (D \rte_\alpha P,j,w) \lori (A,\pi,\rho).
\egloz
\edefin

In other words, the crossed product associated to $(D,P,\alpha)$ is the initial object in the category of covariant representations and their morphisms over $(D,P,\alpha)$. 

The crossed product always exists, but it can be trivial in bad cases. We do not touch upon the question how to detect nontriviality, but remark that there are nice cases treated in the literature where concrete conditions for nontrivial crossed products can be formulated (see for instance \cite{La}).

\bremark
It follows from the definitions that there is a one-to-one correspondence between unital homomorphisms defined on $D \rte_\alpha P$ and covariant representations of $(D,P,\alpha)$ given by 
\bglnoz
  (\Phi: D \rte_\alpha P \lori A) &\loma& (A,\Phi \circ j,\Phi \circ w) \\
  \Phi_{(A,\pi,\rho)} &\lomafr& (A,\pi,\rho).
\eglnoz
Thus, a unital homomorphism $\Phi: D \rte_\alpha P \lori A$ is uniquely determined on $\img(j) \cup \img(w)$ as $\Phi$ can be reconstructed via $\Phi = \Phi_{(A,\Phi \circ j,\Phi \circ w)}$.
\eremark

\bremark
\label{equi}
Given two C*-dynamical semisystems $(D,P,\alpha)$ and $(D',P,\alpha')$ with a $P$-equivariant homomorphism $\varphi: D \lori D'$, every covariant representation $(A,\pi',\rho)$ of $(D',P,\alpha')$ gives rise to a covariant representation $(A,\pi' \circ \varphi,\rho)$ of $(D,P,\alpha)$. This follows from
\bgloz
  \rho(p) (\pi' \circ \varphi(d)) \rho(p)^* = \pi'(\alpha'(p) \circ \varphi(d)) = \pi'(\varphi \circ \alpha(p)(d)) = (\pi' \circ \varphi)(\alpha(p)(d)).
\egloz
Therefore, there exists a uniquely determined morphism of covariant representations, $\Phi: (D \rte_\alpha P,j,w) \lori (A,\pi' \circ \varphi,\rho)$, or, in other words, a homomorphism $\Phi: D \rte_\alpha P \lori A$ such that $\Phi \circ j = \pi' \circ \varphi$ and $\Phi \circ w = \rho$.

If $\varphi$ is an isomorphism, $\varphi^{-1}$ is also equivariant and the crossed products $(D \rte_\alpha P,j,w)$ and $(D' \rte_{\alpha'} P,j',w')$ are isomorphic via the uniquely determined homomorphism $\Psi: D \rte_\alpha P \lori D' \rte_{\alpha'} P$ with $\Psi \circ j = j' \circ \varphi$ and $\Psi \circ w = w'$.
\eremark

\subsection{Semigroup C*-algebras}

From now on, let us fix a left-cancellative semigroup $P$. 

The reduced semigroup C*-algebra can be defined in complete analogy to the group case: For any $p \in P$, consider the operator $\ti{V}_p$ on $\ell^2(P)$ defined by $\ti{V}_p \xi_q = \xi_{pq}$ for all $q \in P$. Since $P$ is left-cancellative, $\ti{V}_p$ is an isometry for any $p \in P$. 

\bdefin
The reduced C*-algebra of $P$ is given by 
\bgloz
  C_r^*[P] \defeq C^*(\menge{\ti{V}_p}{p \in P}) \subseteq \cL(\ell^2(P)).
\egloz
\edefin

Now, the question is how to define the full C*-algebra $C^*[P]$ of $P$. It should be given by universal generators and relations, and we certainly want that $C^*[P]$ is generated by isometries $\menge{\ti{V}_p}{p \in P}$ with $\ti{V}_p \ti{V}_q = \ti{V}_{pq}$ for all $p, q \in P$. But these relations do not contain any information on the range projections of these isometries. Thus, following the idea of Section \ref{R}, we take the smallest family $\cC$ of subsets of $P$ containing $\emptyset$, $P$ and closed under finite intersections as well as left translations (the maps $q \ma pq$ for all $p \in P$). 

\bdefin
The full semigroup C*-algebra associated to $P$ is the universal C*-algebra generated by 
\bgloz
  \text{projections } \menge{\e{X}}{X \in \cC} \text{ and isometries } \menge{\ti{v}_p}{p \in P}
\egloz
such that the following relations hold:
\bgloz
  I.(i) \text{ } \e{\emptyset}=0 \text{, } I.(ii) \text{ } \e{R}=1 \text{, } I.(iii) \text{ } \e{X} \cdot \e{Y} = \e{X \cap Y}
\egloz
\bgloz
  II.(i) \text{ } \ti{v}_p \ti{v}_q = \ti{v}_{pq} \text{, } II.(ii) \text{ } \ti{v}_p \e{X} \ti{v}_p^* = \e{pX}.
\egloz
\edefin

We have $\cC = \menge{(p_1 P) \cap \dotsb \cap (p_n P)}{p_i \in P}$ and thus, $C^*[P]$ is generated by the isometries $\menge{\ti{v}_p}{p \in P}$. Moreover, we certainly have a canonical homomorphism $\pi: C^*[P] \lori C_r^*[P]$ sending generator to generator. Thus, $C^*[P]$ is not trivial.

As a last point, we remark that $D[P] \defeq C^*(\menge{\e{X}}{X \in \cC}) \subseteq C^*[P]$ is a commutative C*-subalgebra and that we can define $C^*[P]$ alternatively as 
\bgloz
  C^*[P] = D[P] \rte_{\alpha_{lt}} P \text{ (in the sense of Definition \ref{defcropro})} 
\egloz
where $P$ acts on $D$ via $\alpha_{lt}(p) \e{X} = \e{pX}$. 

Actually, we can generalize this construction to any family $\cC$ containing $\emptyset$, $P$ and closed under finite intersections as well as left translations. Then, we can consider the family generated by a set $\cF$ of right ideals of $P$, $\cC = \cC(\cF)$. The corresponding C*-algebras will be denoted by $C^*_r[P; \cC]$, $C^*_{r, \cF}[P] = C^*_r[P; \cC(\cF)]$ or $C^*[P; \cC]$, $C^*_{\cF}[P] = C^*[P; \cC(\cF)]$ respectively. These more general constructions will be used to relate this notion to our ring C*-algebras.

Now, let us try to justify our notion of semigroup C*-algebras, with the help of the following three observations:

\subsubsection{Relationship to Nica's construction}

Our construction generalizes Nica's work on quasi-lattice ordered semigroups (see \cite{Ni} or \cite{LaRae}). Namely, it becomes clear that the following property of quasi-lattice ordered semigroups is crucial for Nica's construction:

For any $p_1$, $p_2$ in a quasi-lattice ordered semigroup, $(p_1 P) \cap (p_2 P)$ is either empty or of the form $pP$ for some $p \in P$. 

If this condition holds true in $P$, then one can define the full semigroup C*-algebra by $C^*[P] \defeq D \rte_{\alpha_{lt}} P$ with $D = C^*(\menge{\e{pP}}{p \in P}$. $D$ is multiplicatively closed as the product $\e{p_1 P} \cdot \e{p_2 P} = \e{(p_1 P) \cap (p_2 P)}$ lies in $D$ by the crucial property (compare \cite{Ni}, \cite{LaRae}). 

This definition coincides with ours for quasi-lattice ordered semigroups. The reason it that for these special semigroups, the family $\cC$ will simply be the family $\menge{pP}{p \in P}$. But if we allow more general families of subsets (not just principal right ideals), it is possible to overcome this restriction and to generalize this construction to arbitrary left-cancellative semigroups.

\subsubsection{Semigroup crossed products by automorphisms}

The idea behind our notion of semigroup C*-algebras can be used to clarify the relationship between semigroup crossed products by endomorphisms (see \cite{Mur1}) and automorphisms (see \cite{Mur2} and \cite{Mur3}; we will modify Murphy's definition slightly). 

As we have noted, one serious drawback of semigroup crossed products by endomorphisms as in Section \ref{semicropro} is that we do not have an explicit nontrivial representation of the crossed product at hand. But there is a different notion of semigroup crossed products, due to Murphy:

This time, we look at a unital C*-algebra $D$ and a left-cancellative semigroup $P$ together with a semigrouphomomorphism $\alpha: P \ri \Aut(D)$ (note that $\alpha$ maps to the automorphisms). A covariant representation of $(D,P,\alpha)$ is a triple $(A,\pi,\rho)$ consisting of a unital C*-algebra $A$, a unital homomorphism $ \pi: D \ri A$ and a homomorphism of semigroups $\rho: P \ri \Iso(A)$ with:
\bgloz
  \pi(\alpha_p(d)) \rho(p) = \rho(p) \pi(d) \fa p \in P, d \in D.
\egloz
A morphism $\Phi: (A_1, \pi_1, \rho_1) \lori (A_2, \pi_2, \rho_2)$ of two covariant representations $(A_i, \pi_i, \rho_i)$ ($i = 1, 2$) of $(D,P,\alpha)$ is a homomorphism $\Phi: A_1 \overset{\Phi}{\lori} A_2$ such that $\Phi \circ \pi_1 = \pi_2$ and $\Phi \circ \rho_1 = \rho_2$.

Again, the crossed product $(D \rtimes_{\alpha} P,i,v)$ is defined as the initial object in the category of covariant representations of $(D,P,\alpha)$. But this time, we have a canonical nontrivial representation of $D \rtimes_\alpha P$ in analogy to the left regular representation:

Let $D$ be faithfully represented on a Hilbert space $\cH$. Then, define the following operators on $\ell^2(P,\cH) \cong \ell^2(P) \otimes \cH$: 
\bglnoz
  && \pi(d)(\xi_p \otimes \eta) = \xi_p \otimes \alpha_p^{-1}(d) \eta \\
  && \rho(p)(\xi_q \otimes \eta) = \xi_{pq} \otimes \eta.
\eglnoz
One can check that the covariance relation is fulfilled, so that this gives the desired nontrivial representation of $D \rtimes_\alpha P$. 

So far, this was the definition of Murphy. Now, we will slightly modify this definition. Let $(D \rtm_\alpha P,\ti{i},\ti{v})$ have the universal property as $(D \rtimes_{\alpha} P,i,v)$, but in addition, we ask for one additional property, namely the existence of a homomorphism 
\bgloz
  C^*[P] \lori D \rtm_\alpha P; \ti{v}_p \ma \ti{v}_p.
\egloz
Under this very natural assumption, we get the following relationship: 

\blemma
$D \rtm_\alpha P \cong (D \otimes D[P]) \rte_{\ti{\alpha}} P$ where $\ti{\alpha} \defeq \alpha \otimes \alpha_{lt}$.
\elemma

\bproof
We can use the universal properties to construct mutually inverse homomorphisms. The point is that by our assumption, we can map $D$ and $D[P]$ into $D \rtm_\alpha P$ by homomorphisms with commuting ranges. 
\eproof

\subsubsection{Cuntz algebras and ring C*-algebras as quotients}

The notion of semigroup C*-algebras seems to appear at various places and thereby reveals a unifying character. For instance, as pointed out in \cite{Ni}, the Cuntz algebra $\cO_n$ is a natural quotient of $C^*[\underbrace{\Nz_{0} * \dotsb * \Nz_{0}}_n]$. 

And the ring C*-algebra $\fA_{\cF}[R]$ is - in a canonical way - a quotient of $C^*_{\ti{\cF}}[P_R]$, where $P_R = R \rtimes R \reg$ and $\ti{F} = \menge{I \times (I \cap R \reg) \subseteq R \rtimes R \reg}{I \in \cF}$. Namely, the universal property of $C^*_{\ti{\cF}}[P_R]$ yields a surjection $C^*_{\ti{\cF}}[P_R] \onto \fA_{\cF}[R]$ with $\ti{v}_{(a,b)} \ma u^a s_b$ and $\e{\bigcap_i (a_i, b_i) \cdot (I_i \times (I_i \cap R \reg))} \ma \e{\bigcap_i (a_i + b_i \cdot I_i)}$.

\end{document}